%% file: main.tex
\documentclass[11pt, a4paper]{article}

\input{0_preamble}

\title{\textbf{Distributionally Robust Resource Planning Under Binomial Demand
Intakes}}

\author{ Ben Black\footnote{STOR-i Centre for Doctoral Training, Lancaster
    University, United Kingdom. Email:
    \href{mailto:b.black1@lancaster.ac.uk}{b.black1@lancaster.ac.uk}}
    \footnote{Corresponding author.}\ , Russell Ainslie\footnote{Applied
    Research, BT Technology, Adastral Park, Ipswich, United Kingdom. Email:
    \href{mailto:russell.ainslie@bt.com}{russell.ainslie@bt.com}}\ , Trivikram
    Dokka\footnote{Management Section, Queens Management School, Queens
    University Belfast, United Kingdom. Email:
    \href{mailto:T.Dokka@qub.ac.uk}{T.Dokka@qub.ac.uk}}\ , Christopher
    Kirkbride\footnote{Department of Management Science, Lancaster University
    Management School, United Kingdom. Email:
    \href{mailto:c.kirkbride@lancaster.ac.uk}{c.kirkbride@lancaster.ac.uk}.} }

\begin{document}
%\nocite{*}
\setcitestyle{nosort}

%\pagenumbering{roman} \tableofcontents \newpage

\maketitle
\begin{abstract}
    In this paper, we consider a distributionally robust resource planning model
    inspired by a real-world service industry problem. In this problem, there is
    a mixture of known demand and uncertain future demand. Prior to having full
    knowledge of the demand, we must decide upon how many jobs we \rgf{plan to}
    complete on each day \rgf{in} the plan\rgf{ning horizon}. Any jobs that are
    not completed by the end of their due date incur a cost and become due the
    following day. We present two distributionally robust optimisation (DRO)
    models for this problem. The first is a non-parametric model with a
    phi-divergence based ambiguity set. The second is a parametric model, where
    we treat the number of uncertain jobs due on each day as a binomial random
    variable with an unknown success probability. We reformulate the parametric
    model as a mixed integer program and find that it scales poorly with the
    \rgf{sizes of the} ambiguity and uncertainty sets. Hence, we make use of
    theoretical properties of the binomial distribution to derive fast
    heuristics based on dimension reduction. One is based on cutting surface
    algorithms commonly seen in the DRO literature. The other operates on a
    small subset of the uncertainty set for the future demand. We perform
    extensive computational experiments to establish the performance of our
    algorithms. \rgf{We compare} decisions from the parametric and
    non\rgf{-}parametric models, to assess the benefit of including the binomial
    information.
\end{abstract}

\normalsize
\pagenumbering{arabic}
\input{1_introduction}
\input{2_literature_review}
\input{3_planning_model}
\input{4_experimental_design}
\input{5_results}
\input{6_conclusions}

\Large{\textbf{Acknowledgements}}

\normalsize
We would like to acknowledge the support of the Engineering and Physical
Sciences Research Council funded (EP/L015692/1) STOR-i Centre for Doctoral
Training. We would like to thank BT for their funding, and Mathias Kern and
Gilbert Owusu from BT for their support. We would also like to thank the 4
anonymous EJOR reviewers for their useful comments that have helped us improve
this paper. In addition, we would like to thank Dick den Hertog for his help in
developing the CQP reformulation of the non-parametric version of our model.

\normalsize
\bibliographystyle{apalike}
\bibliography{bib}

\newpage
\input{7_appendices}

\end{document}

%% file: 0_preamble.tex
\usepackage[comma, authoryear, numbers]{natbib}
\usepackage{amsmath,amsfonts,amsthm,bm}
\usepackage[parfill]{parskip}
\usepackage{hyperref}
\usepackage{graphicx}
\usepackage[a4paper, margin=25mm]{geometry}
\usepackage{bm}
\usepackage{bbm}
\usepackage{breqn}
\usepackage[toc,page]{appendix}
\usepackage{tabularx}
\allowdisplaybreaks
\usepackage{multirow}
\usepackage{longtable}
\usepackage{booktabs}
\usepackage{bibentry}
\usepackage[multiple]{footmisc}
\usepackage{algpseudocode}
\usepackage{algorithmicx}
\usepackage{algorithm}
\usepackage{subcaption}
\usepackage{etoolbox}
\usepackage{dsfont}
\usepackage{breqn}
\usepackage{tabulary}
\date{\today}
\graphicspath{{images_tables/}}

\usepackage{amsthm}
\usepackage[dvipsnames]{xcolor}

\usepackage{tocbibind}
\usepackage{comment}
\DeclareMathOperator*{\argmax}{argmax}

\newcommand{\fa}[0]{\ \forall \ }

\newcommand{\E}[0]{\mathbb{E}}

\renewcommand{\P}[0]{\mathbb{P}}
\newcommand{\m}[1]{\mathcal{#1}}

\newcommand{\rgf}[1]{\textcolor{black}{#1}}
\usepackage{mathtools}

\renewcommand{\l}{\left}
\renewcommand{\r}{\right}

%% file: 1_introduction.tex
\section{Introduction}

In this paper, we consider a resource planning problem motivated by a real-world
telecommunications service company. This real problem consists of optimising the
use of a large workforce of service engineers, in the face of a mixture of known
and uncertain jobs.  

\subsection{Problem Setting}\label{sec:problem_setting}

The planning process for a service company is subject to three stages, named
\textit{the three stages of planning}. Each serves a different purpose, covers a
different time horizon, and creates results that feed into the next. The three
stages are \textit{strategic}, \textit{tactical}, and \textit{operational}
planning. Strategic planning covers a period of multiple years, and concerns
long term decisions such as how many employees to be hired and in which skills
they should be trained. Tactical planning concerns a period of weeks or months.
It involves aggregate decisions such as deciding upon the capacity needed in
each period, or how many jobs can and cannot be completed in each period.
Operational planning concerns short-term decisions such as scheduling the
day-to-day activities of the workforce at the individual level. We focus on the
tactical planning stage in this paper. The decisions that we make are at the
aggregate level, i.e.\ we do not plan the specific activities of every worker
but we instead aggregate their availability into a daily capacity value. We are
tasked with planning the use of this capacity to maximise job completions, or
equivalently minimise the number of jobs left incomplete. Since it is typically
not possible to move capacity between days, planners manipulate demand to make
the best use of what they have.

In the telecommunications industry, jobs can be divided into two categories:
repair jobs and installation jobs. Repair jobs correspond to service engineers
being tasked with fixing broken equipment for existing customers, such as
broadband routers and telephone systems. Installation jobs correspond to
engineers installing equipment in order to obtain new customers. For example,
this may be installing new cabling cabinets and networks in order to provide
broadband to a new geographical area. Repair jobs are treated as emergency jobs
and they are given a high priority for completion. Installation jobs are treated
as additional jobs that a company can plan to complete in order to generate more
profit. In this paper, we will consider planning the activities of a
telecommunications workforce carrying out repair jobs. Since breakages in
equipment and services are not planned, these jobs offer a source of
uncertainty. In particular, for any given planning period we have knowledge of a
fixed number of repair jobs that are already in the system at the time of
planning (\textit{workstack jobs}). However, the number of breakages between the
time of planning and the date concerned is subject to uncertainty. The jobs
generated by these future breakages are referred to as \textit{intake jobs}
\rgf{or simply \textit{intakes}}. 

At the time of planning, we have an aggregate capacity value that gives the
number of jobs that our workforce can complete, for each day in a planning
horizon of fixed length. This is obtained from the number of engineers working
on each day, and the number of hours that they will work. By default, we will
use all available capacity on each day to complete jobs that are due on that
day. Furthermore, workstack jobs can be completed on or before their due date,
and completing them early is referred to as \textit{pulling forward}.  However,
the same does not apply to intake jobs. Since the day that they will arrive in
the system is unknown, allowing them to be pulled forward could suggest that
they will be completed before they even arrive. Hence, intake jobs cannot be
pulled forward. If any jobs are still incomplete by the end of their due date,
then they will not leave the system but incur a cost, and become due on the
following day. This is referred to as \textit{rollover}. In this paper, since
capacity is fixed, our model will optimise the pulling forward decision in order
to minimise the total rollover cost over the planning horizon. Pulling forward
can be utilised to free up capacity on due dates that we expect to have high
intake. This helps to reduce rollover and utilise spare capacity.

In the literature on service industry planning models that are closest to ours,
demand uncertainty often results in intractable models due to poor scalability.
Examples of this come from~\cite{AinslieSkills}
and~\cite{Ainslie2018TacticalPO}. In these papers, models had to be solved
heuristically due to their size, even though they were deterministic. However,
the demand uncertainty is still acknowledged. In fact, in some cases the plan is
passed through a predictive model in order to better assess its
performance~\citep{AinslieNN}. The closest model to ours that does model
uncertain demand comes from~\cite{rossthesis}, who used two-stage stochastic
programming models for service industry workforce planning. However, this
methodology requires the assumption that the demand distribution is known, and
this is not an assumption that is reasonable here. The framework that we use to
model our problem is \rgf{d}istributionally \rgf{r}obust \rgf{o}ptimisation
(DRO). This framework allows us to include distributional information in our
models, without full knowledge of the distributions themselves.  

More specifically, we model intakes as binomial random variables where each
distribution is ambiguous. Furthermore, we assume that the intake random
variables for any two days in the plan\rgf{ning horizon} are independent of one
another. We assume that we have access to a forecasting model or expert
knowledge that gives a point estimate of intake and a range of potential values.
Hence, for each day, the number of trials is fixed at the maximum intake.
Therefore, the success probability is the only unknown parameter for each
distribution. This parameter can be estimated through maximum likelihood
estimation, with access to past intake data. Our decision to use the binomial
distribution can be justified by the following three reasons:
\begin{enumerate}
    \item The number of intake jobs due on each day is a discrete quantity and
    any two jobs arriving on the same day arrive independently of one another.
    \item There is a fixed and finite set of values that each intake random
    variable can take. Other discrete distributions such as the Poisson
    distribution are unbounded, and hence not fitting for these random
    variables. 
    \item Apart from naturality, it gives a concise way of modelling the
    uncertainty. We can represent each distribution uniquely by one choice of
    $p$, which is a vector of dimension equal to the number of periods in the
    plan. Using a non-parametric approach would mean having to analyse the
    entire distribution, which is a larger vector that has one entry for every
    realisation of intake.
\end{enumerate}
We emphasize that the binomial assumption is in contrast with much of the DRO
literature, in which distributions are usually non-parametric. The reason for
this is that parametric distributions often lead to intractable models. However,
in the context of our problem, we show that it is possible to derive algorithms
which are both tractable and near-optimal. Binomial and negative binomial
distributions have often been used for demand modelling, particularly in
inventory management. Examples of this include~\cite{collins2004behavior} for a
risk-minimising newsvendor, \cite{gallego2007inventory} for inventory planning
under highly uncertain demand,~\cite{dolgui2008performance} for forecasting
demand in slow-moving inventory systems, and \cite{rossi2014confidence} for
confidence-based newsvendor problems. 

In this paper, we will use the fact that every distribution in the ambiguity set
is binomial in order to find the worst-case expected cost for a fixed pulling
forward decision. In general, our methodology consists of three key steps.
Firstly, we construct a discrete ambiguity set for the parameters of the true
distribution. Secondly, we create a tractable reformulation of the model by
replacing the inner objective with a finite number of constraints. In
particular, there is one constraint for each distribution in the ambiguity set.
Thirdly, we study the objective function as a function of the distribution's
parameters in order to construct a set of extreme \rgf{parameters}. For discrete
distributions, the constraints representing the inner objective will always be
linear. For continuous distributions, this is not necessarily the case. In such
situations, for the second step, one would have to use a linear or quadratic
approximation of the objective function. For example, this could be done using
piecewise linear approximations or sample average approximations. Doing so would
then allow our methodology to be applied.

\subsection{Our Contributions}

We consider a DRO model for a resource planning problem with an unknown number
of intake jobs on each day. Using the problem structure, we model intakes as
binomial random variables and study the resulting DRO model. Due to the use of
the binomial distribution, the problem is considerably harder from a
computational point of view. Our contributions in the paper include the
following:
\begin{enumerate}
    \item  A new framework for solving DRO problems with ambiguity sets
    containing only distributions in the same parametric family as the nominal
    distribution. A comparison of this framework with a common, non-parametric
    framework based on the use of $\phi$-divergences.
    \item Three solver-based algorithms for the parametric model: an optimal and
    a heuristic \rgf{c}utting \rgf{s}urface algorithm \rgf{(named CS\_opt and
    CS, respectively)}, and another heuristic algorithm named Approximate
    Objective (AO) (see Sections~\ref{sec:CS_alg} and~\ref{sec:AO_alg}).
    \rgf{CS}, while not exact, considerably simplifies the main bottleneck step
    of \rgf{CS\_opt}: finding the worst-case distribution for a fixed pulling
    forward decision (referred to as \textit{the distribution separation
    problem}). This makes it much more scalable with the size of the ambiguity
    set.
    \item Extensive computational experiments on a variety of constructed
    instances which show the efficacy of our methods. See
    Section~\ref{sec:results} for these results.
\end{enumerate}

%% file: 2_literature_review.tex
\section{Literature Review}

In this section, we review relevant literature relating to our problem and
problems of a similar nature. In Section~\ref{sec:WP_LR}, we review the
workforce planning literature and highlight the methodologies used there. In
Section~\ref{sec:DRO_LR}, we summarise the recent DRO literature and discuss how
our research differs from it.

\subsection{Workforce and Resource Planning}\label{sec:WP_LR}

Workforce planning models of various forms have been studied in the OR
literature since the mid 1950's, with early papers focussing on creating
tractable deterministic models~\citep{nonlin_55, LP_1960}. Demand uncertainty
has always been discussed in these early papers, with some authors extending
previous models to minimise expected cost rather than cost~\citep{Fetter62}. In
more recent literature, the modelling of uncertain demand has been developed
further. The most common method in the literature has been two-stage stochastic
programming. This methodology was applied to nurse scheduling
~\cite{three_stage_73} and recruitment for a military
organisation~\cite{Price78} in the early literature. More recent examples of
stochastic programming in workforce planning include planning a cyber branch of
the US army~\cite{BASTIAN2020}, and service industry workforce
planning~\cite{Zhu09, rossthesis}. These authors use stochastic programming due
to their assumption that the distribution of the uncertain parameters is known.
When this is not the case, or if the planner is risk-averse, \rgf{r}obust
\rgf{o}ptimisation (RO) can be used to represent demand uncertainty. This
methodology has been used, for example, in healthcare~\citep{HOLTE2013551} and
air traffic control~\citep{RSG_Hertog}.  

Recently, there have also been some applications of DRO to workforce and
resource planning. \cite{LiaoDRO} used DRO for staffing a workforce to take
calls arriving at a call centre at an uncertain rate. The reason for using DRO
was cited as being that the true arrival rates of calls are usually subject to
fluctuations, meaning that the typical stochastic model with a fixed Poisson
distribution was not appropriate. They simulated the DRO solution and the
stochastic programming solution and found that the two had similar costs.
However, the stochastic programming solution violated more model constraints.
\cite{ChenHospitals} also used DRO for workforce planning in a hospital
environment. In particular, they used DRO to determine bed requirements in order
to appropriately manage admissions to the hospital. They use DRO due to the
difficulty in specifying a distribution to describe patient movements in the
hospital, and find that it performs better than a deterministic approach.  

Our resource planning problem deals with the management of both planned and
unplanned jobs. Similar problems exist in other settings, such as scheduling for
gas pipeline maintenance~\citep{FlexibleResource}, and operating room scheduling
in hospitals~\citep{Samudra16}. Particularly, in operating room planning, the
workstack and intake jobs as defined in our model are similar to elective
(inpatient and outpatient) and non-elective (emergency) surgeries. Similarly, in
gas pipeline maintenance the workstack and intake jobs correspond to planned
maintenance jobs and emergency gas leak repairs, respectively. The main
difference between our research and these papers is the choice of performance
measure. For example, \cite{FlexibleResource} use overtime hours as a
performance measure under the assumption that jobs have individual completion
times. However, since our model is for tactical and not operational planning,
jobs and capacity are aggregated. The duration of each job is not modelled
directly. Hence, in our case, the amount of overtime would be inferred by the
number of jobs that could not be completed, i.e.\ rollover. As discussed
by~\cite{Samudra16}, metrics chosen for optimisation differ based on the
underlying context and the stakeholders involved. They emphasize that
traditional metrics such as makespan do not work in presence of both planned and
emergency demands. In our application, the time taken to complete jobs is not of
particular concern. However, leaving jobs incomplete is very costly due to its
effects on customer satisfaction. In industries like telecommunications,
customer satisfaction is of great importance, and hence rollover may be the most
appropriate performance measure.

The literature reviewed here shows that the modelling of uncertain demand in
resource and workforce planning has been the subject of a breadth of research in
the past. It suggests that the most common approach is to employ two-stage
stochastic programming models. However, the assumption that the distribution of
demand is known is not reasonable in our setting. In fact, we only have access
to \rgf{samples of demand and a range of potential values that it can take}.
We do assume, however, that we can take samples of intake in order to estimate
the parameters of its distribution. In addition, there are no recourse actions
in our problem. In such settings, RO and DRO are the only potential solution
approaches. For our problem, a robust model will be shown to lead to more
conservative decisions and large costs. We show this in Appendix~B.2. Hence, we
present a DRO model for our problem, which will extend the previous stochastic
programming approaches to the case where the distribution is not known exactly.
We find that the model is large and complex, due to the size of the sets of
intakes and distributions. Hence, we develop heuristics that apply dimension
reduction to these sets in order to reduce solution times. One algorithm
considers only a small subset of distributions, and the other operates on a
small subset of intakes. While these algorithms perform well on average, they do
sacrifice optimality for speed in some large instances.

\subsection{Distributionally Robust Optimisation}\label{sec:DRO_LR}

DRO combines concepts from robust optimisation and stochastic programming in
order to protect the decision maker from distributional ambiguity. DRO models
are constructed using only limited information on the true distribution of the
uncertain parameters. This information is encoded in an ambiguity set, in which
the true distribution should lie. The earliest type of ambiguity set in the
literature is the moment-based ambiguity set. This set contains all
distributions whose moments satisfy a given set of constraints. The simplest
moment-based sets consider moments to be fixed and known. The moments concerned
have often been the mean and variance. This case was studied
by~\cite{ScarfMinmax} for a newsvendor model. Other papers included models where
the first $m$ moments were known~\citep{Shapiro02}. Authors have also developed
models that did not assume that these values were fixed but that they were known
to lie in an interval or that ordinal relationships between probabilities were
known~\citep{Breton1995}. Other examples of this come from~\cite{Ghaoui03}
and~\cite{LOTFI2018556}, who study a CVaR model where the first two moments are
only known to belong \rgf{to} polytopic or interval sets. Methodologies for
solving moment-based ambiguity \rgf{set} models include reformulation via
bounding the objective function~\citep{ScarfMinmax}, reformulating as a second
order conic program~\citep{Ghaoui03, LOTFI2018556}, sample average
approximations~\citep{Shapiro02}  and sub-gradient
decomposition~\citep{Breton1995}.

The second common methodology for constructing ambiguity sets is using distance
measures. A distance-based ambiguity set contains all distributions that lie
within some pre-prescribed distance of a nominal one. In the literature, many
ways to measure this distance have been studied. For example, many papers have
used the Wasserstein distance. This distance can lead to tractable
reformulations as convex programs~\citep{mohajerin2018data}. Due to this, it has
been used in a number of contexts, such as portfolio selection~\citep{Pflug07},
least squares problems~\citep{MehrotraLS} and statistical
learning~\citep{Lee2015ADA, LeeStatLearn}. 

Another common family of distance measures in DRO has been $\phi$-divergences.
This family contains a number of distance measures, such as the $\chi^2$
distance, variation distance and Kullback-Leibler divergence. Such measures
typically lead to second-order conic programming or even linear programming
relaxations via taking the Lagrangian dual of the inner
problem~\citep{BenTalHertog, Love15}. Due to the convenient reformulations they
yield, $\phi$-divergences have been popular in the DRO literature. There have
been numerous examples of $\phi$-divergences being used to reformulate
\rgf{distributionally robust (DR)} chance-constrained programs as
chance-constrained programs~\citep{hu2013ambiguous, Yanikoglu13,
JiangDataDriven}. Another benefit of $\phi$-divergences is that they can be used
to create confidence sets and enforce probabilistic guarantees.
\cite{BenTalHertog} show how to create confidence sets for the true distribution
based on $\phi$-divergences. This is done by taking a \rgf{maximum likelihood
estimate (MLE)} of its parameters and using the resulting distribution as the
nominal distribution. \cite{Duchi16} use DRO models with $\phi$-divergence
ambiguity sets to construct confidence intervals for the optimal values of a
stochastic program with an ambiguous distribution. Their intervals
asymptotically achieve exact coverage. By studying $\phi$-divergence balls
centred around the empirical distribution, \cite{lam2019recovering} shows that
DRO problems can recover the same standard of statistical guarantees as the
central limit theorem. 

In addition to these papers that consider general $\phi$-divergence functions,
the fact that $\phi$-divergences cover a range of distance measures allows
authors to select those that are most appropriate for their models. For example,
\cite{hanasusanto2013robust} used $\chi^2$ divergence ambiguity sets for a
distributionally robust dynamic programming problem. They used the $\chi^2$
divergence, in particular, because it allows the min-max problems in the dynamic
programming recursion to be reformulated as tractable conic programs. They also
chose this divergence because it does not \textit{suppress} scenarios. In other
words, it does not give scenarios zero probability in the worst-case if they
have non-zero probability under the nominal distribution. The Kullback-Leibler
divergence was also extensively studied by~\cite{hu2013kullback}, who used it
for DR chance-constrained problems. They showed that, under this divergence, if
the nominal distribution was a member of the exponential family then so was the
worst-case distribution.

The literature we have reviewed so far concerns models that can be reformulated
and solved exactly, due to their ambiguity sets being constructed using distance
measures or moment constraints. However, there has also been significant
literature studying general DRO models that are not formulated in this way. In
general, DRO models are semi-infinite convex programs (SCPs). They have a
potentially infinite number of constraints induced by those defining the inner
objective value. Typically, iterative algorithms are used to solve SCP models.
For example, \cite{KortanekHoon93} developed a cutting surface (CS) algorithm
for linear SCP problems with differentiable constraints. This algorithm
approximates the infinite set of constraints with a sequence of finite sets of
constraints. Constraints are iteratively added to the current set considered
until stopping criteria are met. The constraint that is most violated by the
current solution is added at each iteration. In the context of DRO, adding a
constraint corresponds to finding a distribution to add to the current ambiguity
set. This is referred to as solving the \textit{distribution separation
problem}. \cite{Pflug07} later applied this algorithm to DRO models for
portfolio selection under general ambiguity sets. In an extension of
\cite{KortanekHoon93}'s algorithm, \cite{mehrotra2014cutting} developed a CS
algorithm for SCP problems that allowed for non-linear cuts, and did not require
differentiable constraints. CS \rgf{algorithms have} since become a common
apprach to solving DRO problems that are computationally expensive and do not
have tractable reformulations. For example, \cite{RahimianScenarios} applied a
CS algorithm to a DRO model using the total variation distance. They state that
the model becomes expensive to solve to optimality when there are a large number
of scenarios. Another example of its use in the literature is given
by~\cite{Lshaped_DRO}, who used a CS algorithm to solve DR knapsack and server
location problems. \cite{LUO2019} also \rgf{used a CS algorithm} to solve DRO
models under the Wasserstein distance.

Our work differs from the cited literature in two key ways. Firstly, we consider
demand distributions \rgf{belonging} to some parametric family, and enforce that
the worst-case distribution also belongs to this family. We show that the
resulting model can be reformulated as a large MIP. This model becomes slow to
solve for large ambiguity and uncertainty sets. This is due to the large amounts
of computation required and the large number of constraints. Hence, secondly, we
present algorithms that make use of the additional distributional information in
order to solve the parametric model. Among these algorithms is an optimal CS
algorithm, that we will show to be fast for small problems, but to scale poorly
with the size of the ambiguity set. We also contribute a heuristic version of
this CS algorithm, that solves the distribution separation problem at each
iteration over a subset containing only the most extreme parameters. We will
show that this allows us to greatly reduce the time taken to solve the
distribution separation problem. We also show how to construct a confidence set
for the worst-case parameter without the use of $\phi$-divergences. In addition,
we develop the non-parametric model and show how to reformulate it as a
second-order conic program. \rgf{We also compare} the results from the
parametric and non-parametric models to assess the benefit of incorporating the
binomial information. 

%% file: 3_planning_model.tex
\section{\texorpdfstring{Planning Model}{}} \label{sec:models}

In this section, we introduce our planning model and discuss the different types
of ambiguity sets that we will consider. In Section~\ref{sec:notation} we
provide a summary of the notation that will be used. Following this, in
Section~\ref{sec:DRO_model}, we provide the DRO model itself under a general
ambiguity set. In the sections following this, we detail the parametric and
non-parametric versions of the model that will be studied in this paper.

\subsection{Notation and Definitions}\label{sec:notation}

We consider a planning horizon of $L$ periods, which are days in our setting.
The days in the plan are denoted by $\tau \in \{1,\dots,L\}$. The inputs for the
model are defined as follows. For each day $\tau$ we have capacity $c_\tau$,
which gives the number of jobs that we can complete on day $\tau$. The workstack
for day $\tau$ is the number of jobs that are currently due on day $\tau$, and
is denoted $D_\tau$. The workstacks are known at the time of planning. The
intake for day $\tau$ is denoted $I_\tau$. This quantity is the number of jobs
that will arrive between the time of planning and the due date $\tau$ and will
be due on day $\tau$. Each $I_\tau$ is a random variable, and its value is not
realised until the end of day $\tau$. In other words, workstack and intake jobs
represent planned and unplanned/emergency jobs in the terminology used in other
problems.

The rollover for day $\tau$ is the number of jobs that are due on day $\tau$ but
are left incomplete at the end of day $\tau$. This quantity is denoted by
$R_\tau$, which is a random variable due to its dependence \rgf{on $I_\tau$}.
Each unit of rollover on day $\tau$ incurs a cost $a_\tau$. The set of
realisations of the random variable $I_\tau$ is denoted by $\m{I}_\tau =
\{0,\dots, i^{\max}_\tau\}$, and a realisation of $I_\tau$ is denoted by
$i_\tau$. We use suppression of the subscript $\tau$ to represent the vectors of
intakes, workstacks and so on. For example, the vector of workstacks is denoted
by $D = (D_1,\dots,D_L)$.  The set of all realisations of the vector $I$ is
denoted by $\m{I}$. We assume that the set $\m{I}$ is the cartesian product of
the marginal sets, i.e. $\m{I} = \m{I}_1 \times \dots \times \m{I}_L$. In the
language of robust optimisation, $\m{I}$ is referred to as an
\textit{uncertainty set} for $I$. For a realisation $i$ of the vector of intakes
$I$, the corresponding realisation of rollover is denoted by $R^i =
(R^i_1,\dots,R^i_L)$. The objective of our problem is to minimise the
\rgf{worst-case expected} rollover cost by pulling forward jobs. Hence, the
decision vector in our problem is the pulling forward variable, which we denote
by $y$. Jobs can be completed no earlier than $K$ periods prior to their due
date. Therefore, we use $y_{\tau_1, \tau_2}$ to denote the number of jobs pulled
forward from period $\tau_1 \in \{2,\dots,L\}$ to period $\tau_2 \in \{\tau_1 -
K,\dots, \tau_1 - 1\}$. This corresponds to completing $y_{\tau_1, \tau_2}$
additional jobs on $\tau_2$ that are due \rgf{on} $\tau_1$.

\subsection{General Distributionally Robust Model}\label{sec:DRO_model}

We now consider the distributionally robust planning model, which is defined as
follows. Denote by $\m{P}$ a general \textit{ambiguity set} of intake
distributions, such that every distribution $P \in \m{P}$ assigns a probability
to every possible intake $i \in \m{I}$. Our model aims to minimise the
worst-case expected rollover cost by selecting the value of $y$. The model is
shown in~\eqref{eq:Lday_obj}-\eqref{eq:Lday_last}.
\begin{align}
    \min_{y, R} & \max_{P \in \m{P}} \sum_{\tau=1}^L a_\tau \E_P(R_\tau)
    \label{eq:Lday_obj}\\
    \text{s.t.}  &\sum_{\tau_2 = \max\{\tau_1 - K, 1\}}^{\tau_1 - 1} y_{\tau_1,
    \tau_2} \le D_{\tau_1} \fa \tau_1 = 2,\dots, L,\label{eq:move_to}\\
    &\sum_{\tau_1 = \tau_2 + 1}^{\min\{\tau_2 + K, L\}} y_{\tau_1, \tau_2} \le
    \max\{c_{\tau_2} - D_{\tau_2}, 0\} \fa \tau_2 = 1,\dots, L - 1,
    \label{eq:move_from}\\
    &R^i_1 \ge i_1 + \sum_{\tau_1 = 2}^{\min\{1 + K, L\}} y_{\tau_1, 1} - \l(c_1
    - D_1\r) \fa i \in \m{I}, \label{eq:R_1}\\
    &R^i_\tau \ge R^i_{\tau - 1} + i_\tau + \sum_{\tau_1 = \tau + 1}^{\min\{\tau
    + K, L\}} y_{\tau_1, \tau} - \l(c_\tau - D_\tau + \sum_{\tau_2 = \max\{\tau
    - K, 1\}}^{\tau -1} y_{\tau, \tau_2}\r) \nonumber\\
    & \fa \tau = 2, \dots,L-1 \fa i \in \m{I},\label{eq:R_tau}\\
    &R^i_L \ge R^i_{L - 1} + i_L - \l(c_L - D_L + \sum_{\tau_2 = \max\{L - K,
    1\}}^{L - 1} y_{\tau, \tau_2}\r) \fa i \in \m{I},\label{eq:R_L}\\
    &y_{\tau_1, \tau_2} \in \mathbb{N}_0 \fa \tau_1, \tau_2,\\
    &R^i_\tau \ge 0 \fa \tau =1,\dots,L \fa i \in \m{I}.\label{eq:Lday_last}
\end{align}
The general idea in calculating rollover in the $L$-day model is as follows. For
a given day $\tau$, we first compute the number of jobs to be completed on day
$\tau$. To compute this, we take the rollover from day $\tau - 1$ and day
$\tau$'s intake as a baseline number of jobs. Then we add the number of jobs
pulled forward to day $\tau$, i.e.\ $\sum_{\tau_1 = \tau + 1}^{\min\{\tau + K,
L\}} y_{\tau_1, \tau}$. We then compute the capacity that can be used to
complete these jobs. This is done by taking the capacity $c_\tau$ and
subtracting the capacity required to complete those workstack jobs that are not
pulled forward from day $\tau$, i.e.\ $D_\tau - \sum_{\tau_2 = \max\{\tau - K,
1\}}^{\tau -1} y_{\tau, \tau_2}$. If the remaining capacity is enough to
complete all jobs on $\tau$, then the rollover is zero. Otherwise, the rollover
is the number of jobs left incomplete.  

Constraints~\eqref{eq:move_to} and~\eqref{eq:move_from} provide upper bounds on
the pulling forward totals. Constraint~\eqref{eq:move_to} ensures that no jobs
are pulled forward if they cannot be completed on the day to which they are
moved. Constraint~\eqref{eq:move_from} ensures that only workstack jobs can be
pulled forward, and that a job cannot be pulled forward multiple times in order
to be pulled forward more than $K$ days. Constraint~\eqref{eq:R_1} reflects that
jobs cannot be pulled forward from day 1 and hence we only subtract those jobs
pulled forward to day 1 from its remaining capacity. We do not reduce rollover
by pulling forward from it. Similarly, constraint~\eqref{eq:R_L} reflects that
jobs cannot be pulled forward to the final day of the plan. Hence, we only pull
forward from this day and not to this day. For every other day,
constraint~\eqref{eq:R_tau} captures that we can pull forward to \textit{and}
from said day. We therefore add \textit{and} subtract jobs from its capacity to
calculate the rollover. 

\subsection{Non-parametric DRO Model} \label{sec:NP_model}

The non-parametric model is defined by ambiguity sets $\m{P}$ containing
distributions $P$ that are not necessarily parametric. To be specific, $\m{P}$
can be any subset of the set of all distributions over the set of intakes, i.e.\
$\m{P} \subseteq \l\{P \in [0,1]^{|\m{I}|}: \sum_{j=1}^{|\m{I}|} P_j = 1\r\}$.

\subsubsection{Phi-divergence Based Ambiguity Sets}\label{sec:phi_sets}

As discussed earlier in the paper, it is common to define $\m{P}$ using
$\phi$-divergences. Adopting similar notation to that of~\cite{Love15}, suppose
that $P$ and $Q$ are two probability distributions. We define a
$\phi$-divergence $d_{\phi}$ for $\phi$-divergence function $\phi$ as:
\begin{equation}
    d_{\phi}(P, Q) = \sum_{j=1}^n Q_j \phi\l(\frac{P_j}{Q_j}\r),
\end{equation}
where $\phi$ is a convex function on the non-negative reals. This function
measures the distance between $P$ and $Q$. \rgf{In what follows, $Q$ will be
treated as a nominal distribution.} Furthermore, we denote by $\phi^*$ the
\textit{conjugate} of $\phi$, which can be found via~\eqref{eq:conj}.
\begin{equation}\label{eq:conj}
    \phi^*(s) = \sup_{t \ge 0}\{st - \phi(t)\} 
\end{equation}
The conjugate will be useful when finding reformulations later in the paper.
Given a nominal distribution $Q$, we can define $\m{P}$ as the set of all
distributions $P$ that lie within some pre-prescribed distance from $Q$ as
measured by the $\phi$-divergence. In other words, we can use:
\begin{equation}\label{eq:NP_conf_set}
     \m{P}_\rho = \l\{P \in [0,1]^{|\m{I}|}: \sum_{j=1}^{|\m{I}|} P_j = 1,
     d_\phi(P, Q) \le \rho\r\}.
\end{equation}
As described by~\cite{BenTalHertog}, this formulation of the ambiguity set
allows us to choose $\rho$ such that $\m{P}$ is a confidence set for the true
distribution. Suppose that the true distribution $P^0$ lies in a parameterised
set $\{P^{\theta}\ | \ \theta \in \Theta\}$, such that the true value of
$\theta$ is $\theta^0$. Also suppose that we take $N$ samples of intake from
$P^0$ and take an MLE $\hat{\theta}$ of $\theta^0$. Then, if we choose $\rho$
using~\eqref{eq:rho}, the set $\m{P}_\rho$ is an approximate $100(1-\alpha)\%$
confidence set for $P^0$ around $\hat{P} = P^{\hat{\theta}}$.
\begin{equation}\label{eq:rho}
    \rho = \frac{\phi''(1)}{2N} \chi^2_{k, 1-\alpha}.
\end{equation}
In~\eqref{eq:rho}, $k$ is the dimension of $\Theta$ and $\chi^2_{k, 1-\alpha}$
is the $100(1-\alpha)^{\text{th}}$ percentile of the $\chi^2$ distribution with
$k$ degrees of freedom. There are many choices for the choice of
$\phi$-divergence function, and some examples can be found in the paper
by~\cite{BenTalHertog}.

\subsubsection{Reformulation with \texorpdfstring{Modified
$\chi^2$-divergence}{}} \label{sec:reformulation}

In our model, we will use the modified $\chi^2$ distance as our
$\phi$-divergence. This uses the $\phi$-divergence function $\phi_{m\chi^2}(t) =
(t-1)^2$ and is defined in~\eqref{eq:mod_chisq}.
\begin{equation}\label{eq:mod_chisq}
    d_{\phi_{m\chi^2}}(P, Q) = \sum_{j=1}^{n}\frac{\l(P_j - Q_j\r)^2}{Q_j}.
\end{equation}
Here, $n$ is the number of potential values of the uncertain parameters. In our
problem, we have $n = |\m{I}|$. We choose this function for the following
reasons. Firstly, it leads to a convex quadratic programming (CQP)
reformulation. Secondly, squared deviations from the nominal distribution are
represented as a proportion of the nominal distribution's value. This means that
small deviations from the nominal distribution can still lead to a large term in
the sum in~\eqref{eq:mod_chisq}. When $n$ is large, most values of $Q_j$ will be
small, and this will help identify significant deviations from small nominal
values. Other choices of $\phi$-divergences that lead to CQP reformulations,
such as the $\chi^2$ distance, Hellinger distance and the Cressie-Read distance,
do not have the normalisation effect given by dividing each term by $Q_j$.
Following~\cite{BenTalHertog}, defining $s_j = \frac{\sum_{\tau=1}^L a_\tau
R^{i^j}_\tau - \nu}{\lambda}$, \rgf{we find} the following CQP reformulation of
our full model \rgf{with $\m{P} = \m{P}_\rho$}:
\begin{align}
    \min_{y, R, \lambda, \nu, z, u} \ & \l\{ \lambda (\rgf{\rho} - 1) + \nu +
    \frac{1}{4}\sum_{j=1}^n Q_j u_j \r\}\label{eq:final_reform_np_obj},\\
    \text{s.t. } &\eqref{eq:move_to}-\eqref{eq:Lday_last},\\
    &\sqrt{4z^2_j + (\lambda - u_j)^2} \le (\lambda + u_j) \fa j = 1, \dots, n\\
    &z_j \ge \sum_{\tau=1}^L a_\tau R^{i^j}_\tau - \nu + 2\lambda \fa j =
    1,\dots,n\\
    &z_j \ge 0 \fa j = 1,\dots, n.\\
    &\lambda \ge 0.
\end{align}
In this formulation, $z_j$ and $u_j$ for $j=1,\dots,n$ are dummy variables
defined to ensure that the model is a CQP model. A full derivation of this
reformulation can be found in Appendix~A, along with how to extract the
worst-case distribution from its solution.

\subsection{Parametric DRO Model}

In this section, we detail a parametric version of the DRO planning model. This
is a new modelling framework for DRO problems that allows the ambiguity set to
contain only distributions that are members of the same parametric family as the
true distribution. This is useful in cases where we know beforehand which family
the true distribution lies in, because it ensures that the worst-case
distribution implied by the model is also in this family.

\subsubsection{Implications of Parametric Ambiguity Sets}\label{sec:param_sets}

Recall from Section~\ref{sec:phi_sets} that we can use $\phi$-divergences to
create confidence sets when we know that the true distribution lies in some
parametric family $\m{P}_{\Theta} = \{P^{\theta}\ | \ \theta \in \Theta\}$. The
resulting confidence set~\eqref{eq:NP_conf_set}, however, does not only contain
distributions in this family. Therefore, there is no guarantee that the
worst-case distribution will lie in this family and hence no guarantee that it
is even a distribution that could be equal to $P^0$. Our methodology involves
explicitly using the set $\m{P}_{\Theta}$ in our DRO model instead, which
eliminates potential worst-case distributions that are not in the same family as
the true distribution. Suppose that we take the ambiguity set given by $\m{P} =
\m{P}_{\Theta}$. 

The methodology in Section~\ref{sec:reformulation} relies on being able to
represent the requirement that $P \in \m{P}$ in the constraints of the model.
However, representing $P \in \m{P}_{\Theta}$ in the constraints is more
challenging. In the case where $\m{P}_{\Theta}$ represents a set of discrete
parametric distributions, e.g.\ binomial or Poisson, the requirement might be
represented by:
\begin{equation}
    P_j = f(i^j\ | \ \theta) \text{ for some } \theta \in \Theta,
\end{equation}
where $f$ is the probability mass function (PMF) of $I$ and $i^j$ is the
$j^{\text{th}}$ realisation of intake. The only reasonable way that one might
attempt to include this in the model is to treat $\theta$ as a dummy variable,
and replace $P_j$ in the objective with $f(i^j \ | \ \theta)$. However, most
PMFs as functions of their parameters are either high order polynomials (such as
binomial) or include exponential functions (such as Poisson). Including them in
the model through the objective function will hence make the model intractable.
As an example, consider our model with independent intakes and $I_\tau \sim
\text{Bin}(i^{\max}_\tau, p_\tau)$ for $\tau=1,\dots,L$. The objective of the
inner problem becomes:
\begin{equation}\label{eq:obj_w_pmf}
    \max_{p \in \Theta} z = \sum_{\tau=1}^L \sum_{i \in \m{I}} a_\tau R^i_\tau
    \prod_{l = 1}^L \binom{i^{\max}_l}{i_l} p^{i_l}_l (1-p_l)^{i^{\max}_l -
    i_r}.
\end{equation}
Treating this as a \rgf{non-linear program}, we might consider solving using the \rgf{Karush-Kuhn-Tucker (KKT)} conditions. The derivative of
the objective function in~\eqref{eq:obj_w_pmf} \rgf{with respect to} $p_k$ is:
\begin{equation}
\sum_{\tau, i} a_\tau R^i_\tau  \binom{i^{\max}_k}{i_k}\l(i_k p_k^{i_k -
1}(1-p_k)^{i^{\max}_k - i_k} - p^{i_k}_k(i^{\max}_k - i_k)(1-p_k)^{i^{\max}_k -
i_k - 1}\r)\prod_{l \neq k}f_l(i_l),
\end{equation}
for each $k \in \{1,\dots,L\}$, where $f_l$ is the PMF of $I_l$. Choosing a
vector $p$ such that $p_\tau < 1$ for all $\tau$ and all derivatives are equal
to zero is a challenging task. This would need to be done numerically, and hence
would not result in a tractable objective function for our outer model.
Furthermore, using a $\phi$-divergence to define $\Theta$ would not result in a
tractable reformulation. This would involve using an ambiguity set for $p$ of
the form:
\begin{equation}
    \Theta = \{p \in [0,1]^L: d_\phi(p, q) \le \rgf{\rho}\},
\end{equation}
where $q$ is the success probability vector corresponding to the nominal
distribution $Q$. Now consider the methodology in
Section~\ref{sec:reformulation}. This methodology relies on the objective
function being separable over $j$ (see Appendix A). Following the same steps but
with the objective in~\eqref{eq:obj_w_pmf}, we arrive at the following dual
objective:
\begin{equation}
    \min_{\lambda \ge 0}\l\{ \lambda_0 \rgf{\rho} + \lambda_1\max_{p \ge
    0}\sum_{\tau=1}^L\l( \sum_{i \in \m{I}} a_\tau R^{i}_\tau \prod_{l = 1}^L
    f_l(i_l)  - \lambda_0 q_\tau \phi\l(\frac{p_\tau}{q_\tau}\r ) + \lambda_1
    (1-p_\tau)\r)\r\}.
\end{equation}
Due to the product over $l$ inside the $\max_{p \ge 0}$ operator (which contains
each success probability), we see that this objective is not separable over
$\tau$. Thus, the remaining steps in creating a tractable reformulation cannot
be carried out. This holds not only for independent distributions, but for any
\rgf{distribution} where the PMF of $I$ depends on more than one $p_\tau$. 

Hence, our methodology is as follows. Instead of treating the parameter $\theta$
as a vector of decision variables, we represent it using a discrete, finite set
of potential values. In other words, we assume that $\Theta$ is a discrete and
finite set. This allows us to represent the distributional ambiguity via a
finite set of constraints that are linear in the rollover variables. The
resulting model has one additional constraint for every $\theta \in {\Theta}$,
but remains a tractable mixed integer program (MIP). We detail the MIP
reformulation of the parametric model in Section~\ref{sec:IP_reformulation}.

\subsubsection{Mixed Integer Programming
Reformulation}\label{sec:IP_reformulation}

To solve this model, we can reformulate it as an MIP as follows. Firstly, we
replace the set $\m{P}_{\Theta}$ with $\Theta$ and optimise over the parameters
$p$ directly. Since there is a one-to-one mapping between $\theta$ and
$P^{\theta}$, the objective becomes:
\begin{equation}\label{eq:para_obj_theta}
    \min_{y, R}\max_{\theta \in \Theta} \sum_{\tau=1}^L a_\tau
    \E_{\theta}\l(R_\tau\r).
\end{equation}
Next, we define a dummy variable $t$ to represent the worst-case expected cost
for a given $y$. Since the set $\Theta$ is a discrete set, we can enforce the
requirement that $t = \max_{\theta \in \Theta} \sum_{\tau=1}^L a_\tau
\E_\theta(R_\tau)$ a set of linear constraints. Hence, the MIP reformulation of
the DRO model is given by:
\begin{align}
    \min_{y, R, t} \ &t \label{eq:IP_obj}\\ \text{s.t. } &
    \eqref{eq:move_to}-\eqref{eq:Lday_last}, \\ &t \ge \sum_{\tau=1}^L a_\tau
    \E_\theta(R_\tau) \fa \theta \in \Theta,\label{eq:IP_t_constraint} 
\end{align}
This model can be very slow to build and solve. This is mostly due to the amount
of computation required to build the model and its constraints. The constraint
for $t$ requires us to compute the distribution $P^\theta$ for every $\theta \in
\Theta$. Due to the sizes of $\Theta$ and $\m{I}$, this can be very slow. To see
this, consider an example with $|\Theta| \ = 3883$ distributions and $|\m{I}| \
= 20000$ potential intakes. Suppose also that the intakes are independent. Then,
for each of 3883 distributions we would need to compute a product of $L$ PMF
values, for each of 20000 intakes. This means computing $L \times 3883 \times
20000 = L \times(77.66\times 10^6)$ PMF values. Furthermore, the model has
$L|\m{I}|$ rollover variables and constraints, and $|\Theta|$ expected value
constraints. This also makes the model slow to build and solve for large
instances. For this instance with $L=5$, this corresponds to 103,878 additional
constraints, when compared with the deterministic model. Our heuristics
therefore employ dimension reduction techniques to make them more tractable. 

\subsection{Binomial Intakes and Ambiguity Sets}

As previously discussed, we will assume that the intakes in our problem our
binomially distributed. In other words, we assume that $I_\tau \sim
\text{Bin}(i^{\max}_\tau, p_\tau)$. We assume that $\m{I}$ is provided to us
prior to model building, either by a prediction model or expert knowledge. The
true set in which we know that the true $p$, denoted $p^0$, must lie is
$[0,1]^L$. As detailed in Section~\ref{sec:param_sets}, we will however use a
finite, discrete subset of $[0,1]^L$ as an ambiguity set for our model. We
consider a discretisation of $[0,1]^L$ of the form given in~\eqref{eq:U^P^0},
where $n_{\text{probs}}$ is chosen by the planner, and details the fineness of
the discretisation.
\begin{equation}\label{eq:U^P^0}
    \Theta_{\text{base}} = \l\{\frac{j}{n_{\text{probs}}} \ \Bigg | \ j =
    0,\dots,n_{\text{probs}}\r\}^L
\end{equation}
Secondly, we assume that we have access to $N$ samples of past intake data, from
which we can take an MLE $\hat{p}$ of $p^0$. The corresponding distribution is
given by $\hat{P}$, which has mean vector $\hat{i} = \hat{p} i^{\max}$. Given
the MLE $\hat{p}$, we consider only $p \in \Theta_{\text{base}}$ that can be
considered \textit{close} to $\hat{p}$. As mentioned earlier, it is common in
the non-parametric DRO literature to use $\phi$-divergences to measure the
distance between two distributions. The main reason for this is that it results
in tractable reformulations via dualising the inner problem. However, since our
approach does not entail dualising the inner problem, this benefit does not
apply to us. Another reason for using $\phi$-divergences is that they allow us
to create confidence sets for the true distribution. However, this is based on
applying the $\phi$-divergence to the distributions themselves, not to the
parameters. We could construct a confidence set for $p^0$ by first constructing
a confidence set for $P^0$ and then creating $\Theta$ by extracting the
parameters of each distribution in the confidence set. However, this would
entail computing the corresponding distribution for every $p \in
\Theta_{\text{base}}$, which is a large computational task. Hence, we do not use
$\phi$-divergences for the parametric model. We can, however, construct a
confidence set for $p^0$ without using $\phi$-divergences and without needing to
compute each distribution $P^p$. Since $\hat{p}$ is an MLE of $p^0$ based on $N$
samples from the true intake distribution,
by~\cite{MillarRussellB2011MLEa}\rgf{,} for large $N$ we have:
\begin{equation}
    (\hat{p}_\tau - p^0_\tau) \sim \mathcal{N}\l(0,\frac{\hat{p}_\tau
    (1-\hat{p}_\tau)}{N i^{\max}_\tau}\r),
\end{equation}
approximately. Therefore, by independence of the $L$ different MLE's, we have
that:
\begin{equation}
    \sum_{\tau=1}^L \frac{N i^{\max}_\tau}{\hat{p}_\tau
    (1-\hat{p}_\tau)}(\hat{p}_\tau - p^0_\tau)^2 \sim \chi^2_L,
\end{equation}
approximately. \rgf{Hence,} we have the following \rgf{discretisation of an}
approximate $100(1-\alpha)\%$ confidence set for $p^0$ around $\hat{p}$:
\begin{equation}\label{eq:conf_set_fast}
    \Theta_\alpha = \l\{p \in \Theta_{\text{base}}: \sum_{\tau=1}^L N
    i^{\max}_\tau\frac{(\hat{p}_\tau - p_\tau)^2}{\hat{p}_\tau (1-\hat{p}_\tau)}
    \le \chi^2_{L, 1-\alpha}\r\}.
\end{equation}
This may yield a different ambiguity set to the one obtained using the
$\phi$-divergence method. This is because they are two different approximations
of the same set. 

\subsection{Solver-based Solution Algorithms}\label{sec:CS_AO}

As described in Section~\ref{sec:IP_reformulation}, the model can be solved to
optimality by reformulating it as a mixed integer program. However, when
$\Theta$ and $\m{I}$ are large, this model has a large number of constraints and
decision variables. This can make it very slow to solve. Hence, we develop three
dimension reduction algorithms in order to reduce the effects of the sizes of
these sets on solution times. In Section~\ref{sec:CS_alg}, we discuss two
cutting surface (CS) algorithms. The first is an optimal CS algorithm that also
scales poorly with the size of $\Theta$. The second is a heuristic CS algorithm
that applies dimension reduction to $\Theta$. Then, in Section~\ref{sec:AO_alg},
we describe our \rgf{final algorithm, Approximate Objective (AO),} that applies
dimension reduction to $\m{I}$.

\subsubsection{Cutting Surface Algorithms} \label{sec:CS_alg}

In this section, we describe our adaptation of the CS algorithm detailed in the
literature review, which has been commonly used in the DRO literature. The
algorithm has a number of different forms, but the one that we base our
adaptation on is that from the review paper by~\cite{rahimian2019}. The general
idea of the algorithm is as follows. In order to deal with the large number of
constraints implied by the ambiguity set, the algorithm uses the following
steps. We start with a singleton set containing one distribution, and solve the
problem \rgf{over} this ambiguity set. Then, for the generated pulling forward
solution, we find the worst-case distribution over the entire ambiguity set. We
add this distribution to the \rgf{current} ambiguity set and then repeat the
above steps. This procedure repeats until optimality criteria are met. 

In more detail, suppose that we have some initial subset $\Theta^{\rgf{1}}
\rgf{= \{p^0\}}$ of our set of distributions $\Theta$ and we solve the full
model with ambiguity set $\Theta^{\rgf{1}}$, to get an optimal decision
$y^{\rgf{1}}$. Then, we find the worst-case parameter, $p^{\rgf{1}} \in \Theta$,
for the fixed solution $y^{\rgf{1}}$, and add it to our set to create a new
subset $\Theta^{\rgf{2}} = \Theta^{\rgf{1}} \cup \{p^{\rgf{1}}\}$ of $\Theta$.
We then solve the model with ambiguity set $\Theta^{\rgf{2}}$, and repeat. We
stop the algorithm when we have reached $\varepsilon$-optimality, i.e.\ if the
solution from the full problem at iteration $k$, $y^k$, gives a worst-case
expected cost over $\Theta^k$ that is within $\varepsilon/2$ of the worst-case
expected cost for $y^k$ over $\Theta$. The algorithm returns an
$\varepsilon$-optimal solution to the DRO model in a finite number of
iterations. The issue with this version of CS is that, even if $y$ is fixed at
$y^k$, finding the true  worst-case distribution $p^k$ can be a very cumbersome
task. To do so, this problem is often treated as an LP. In our case, we can
simply enumerate all distributions in $\Theta$. Even though this is not a
difficult task, it requires a significant amount of computation due to the
necessity of calculating the distributions themselves. 

From now on, we refer to the optimal CS algorithm described above as CS\_opt. We
will show that this algorithm suffers from poor scaling with respect to the size
of $\Theta$. In order to reduce the computational burden, we apply dimension
reduction to $\Theta$. Particularly, we use the simple observation that
$\E_p(I_\tau) = i^{\max}_\tau p_\tau$ is increasing in $p_\tau$ to construct a
set of extreme \rgf{parameters}. Intuitively, this result suggests that a higher
success probability also leads to no-less expected rollover, due to the fact
that $R^i_\tau$ is increasing in $i_\tau$. Hence, we construct a set of
probability vectors such that at least one value is maximised. If this is not
the case, then one value can be increased and this would cause higher expected
rollover for that day. Furthermore, we also assume that the total success
probability is maximised given that one value is maximised. This is to ensure
that we take the most extreme probability vectors over all those such that one
success probability is maximised. Mathematically, we define the set of extreme
parameters as follows. Define $p^{\max}_\tau = \max_{p \in \Theta} p_\tau$ for
$\tau = 1,\dots,L$ and find the set of parameters such that one value is
maximised: 
$$\Theta^{\max}_\tau = \l\{p \in \Theta: p_\tau = p^{\max}_\tau\r\} \text{ for }
\tau = 1,\dots,L.$$ For each $\tau$, construct a set of the most extreme
parameters in $\Theta^{\max}_\tau$ and take the union of these sets to form
$\Theta^{\text{ext}}$:
$$\Theta^{\text{ext}}_\tau = \argmax_{p \in \Theta^{\max}_\tau}\l\{\sum_{k=1}^L
p_k\r\}, \quad \Theta^{\text{ext}} = \bigcup_{\tau = 1}^L
\Theta^{\text{ext}}_\tau.$$

In order to reduce the computation required, our heuristic CS algorithm
\rgf{(referred to as CS)} solves the distribution separation problem over
$\Theta^{\text{ext}}$, rather than the entire ambiguity set $\Theta$. The
general framework for both of our CS algorithms is given below, where CS\_opt
uses $\tilde{\Theta} = \Theta$ in step 2(b)\ and CS uses $\tilde{\Theta} =
\Theta^{\text{ext}}$.
\begin{enumerate}
    \item Compute ambiguity set $\tilde{\Theta}$ and initialise
    $\Theta^{\rgf{1}} = \{p^0\}$, where $p^0 = \hat{p}$ for example.
    \item For $k = \rgf{1},\dots,k^{\max}$:
    \begin{enumerate}
        \item Solve the model to optimality using \rgf{ambiguity set} $\Theta^k$
        to generate solution $(y^k, t^k)$ where $t^k$ is worst-case expected
        cost of $y^k$ over the set $\Theta^k$ passed to the model.
        \item Solve distribution separation problem $\max_{p \in \tilde{\Theta}}
        \sum_{\tau=1}^L a_\tau \E_p(R_\tau\ |\ y = y^k)$ to get solution $p^k$:
        \begin{enumerate}
            \item For $p \in \tilde{\Theta}$\rgf{,} calculate $C_p =
            \sum_{\tau=1}^L a_\tau \E_p(R_\tau \ | \ y = y^k)$.
            \item Choose $p^k$ such that $C_{p^{k}} = \max_{p \in
            \tilde{\Theta}}(C_p)$.
        \end{enumerate}
        \item If $C_{p^k} \le t^k + \frac{\varepsilon}{2}$ or $p^k \in \Theta^k$
        then stop and return solution $(y^k, p^k)$.
        \item Else, set $\Theta^{k+1} = \Theta^k \cup \{p^k\}$ and $k = k + 1$.
    \end{enumerate}
\end{enumerate}
The logic behind 2(c), where we check if $p^k \in \Theta^k$, is that calculation
differences might cause $t^k$ and $C_{p^k}$ to differ by more than
$\frac{\varepsilon}{2}$ when they should be equal. Solvers use some dimension
reduction techniques when building and solving their models. This can lead to
objective values that are not the same as the ones given by the function used in
2(b), even for the same arguments. This stopping criterion is also used in the
CS algorithms by~\cite{Pflug07} and~\cite{Lshaped_DRO}. We now explain why the
condition cannot cause early stopping. Firstly, assume that $\hat{p}$ is not a
worst-case parameter for $y^k$ in $\Theta^k$, i.e.\ it did not give a cost of
$t^k$. Since $p^k$ is generated by the distribution separation problem, it is a
worst-case parameter for $y^k$ over the entire set $\tilde{\Theta}$. If we also
have $p^k \in \Theta^k$ then we have the following two facts. Firstly, we have
$\Theta^k \setminus \{\hat{p}\} \subseteq \tilde{\Theta}$ and so $p^k$ is
necessarily worse than every $p \in \Theta^k \setminus \{\hat{p}\}$. Secondly,
$p^k$ must be worse than $\hat{p}$, because otherwise $\hat{p}$ would be a
worst-case parameter in $\Theta^k$. Hence, $p^k$ is a worst-case parameter in
$\Theta^k$, i.e.\ $C_{p^k} = t^k < t^k + \frac{\varepsilon}{2}$. Now suppose
that $\hat{p}$ is a worst-case parameter in $\Theta^k$. If $p^k \in \Theta^k$
then we must have $C_{p_k} \le t_k < t^k + \frac{\varepsilon}{2}$ since
$\hat{p}$ is worse than $p^k$. Hence, whenever $p^k \in \Theta^k$ occurs, the
first stopping criterion should also be met. 

\subsubsection{Approximate \rgf{O}bjective Algorithm}\label{sec:AO_alg}

The final algorithm that we describe is named Approximate \rgf{O}bjective (AO).
When solving the model to optimality, we are required to compute the
distribution $P^p$ for each $p \in \Theta$. For each intake $i^j \in \m{I}$ we
can easily compute:
\begin{equation}
    \max_{p \in \Theta} P^p_j = \max_{p \in \Theta} \P(I = i^j\ |\ p),
\end{equation} 
and then we can consider a new set of intakes in the model defined by:
\begin{equation}
    \tilde{\m{I}} = \l\{i \in \m{I}: \max_{p \in \Theta} P^p_j > \beta\r\}
\end{equation}
where $\beta$ is our minimum intake probability. By tuning $\beta$, we are
removing intakes from our set that are very unlikely. When solving the model, we
are approximating the expected value by removing some small terms. Since the
intakes removed have low probability, this approximation should be strong. We
simply solve the MIP reformulation with the full set $\Theta$ of
\rgf{parameters} but over the reduced set $\tilde{\m{I}}$ of intakes. For this
paper, we use $\beta = 10^{-3}$ as our initial testing showed that this value
led to good improvements in computation \rgf{time}.

\subsection{Example: A Two-day Problem}

In order to illustrate the logic behind our algorithms, we now give an example
of their use for a two-day version of our model. Since there is only one
feasible pair of days that we can pull forward jobs between, i.e.\ $(2,1)$,
there is now only one decision variable. We refer to this decision variable as
$y = y_{2,1}$. The two-day model is given
by~\eqref{eq:Lday_obj}-\eqref{eq:Lday_last} with $L=2$ and $K=1$.

Suppose that we have $c = (30, 10)$, $D = (5, 20)$, $i^{\max} = (20, 20)$ and $a
= (1,1)$. This gives $|\m{I}| \ = 21^2 = 441$. We construct a 99.5\% confidence
set for $p^0$ using  $\alpha = 0.005$, $N = 10$ and $n_{\text{probs}} = 100$.
This gives $|\Theta| \ = 305$, and we find that the maximum values of $p_1$ and
$p_2$ are both 0.84. This suggests that the above model has $2 \times 441 = 882$
rollover constraints and variables, 81 expected value constraints and 2 pulling
forward constraints. Hence, it has 1189 constraints and 884 decision variables.
We solve this model to optimality in 2.6 seconds, to find the optimal $y$ to be
$y^{\text{P}} = 9$ and the worst-case $p$ to be $p^{\text{P}} = (0.82, 0.82)$
with an expected cost of $z^{\text{P}} = 19.2$. 

When we solve this model with CS, we find that $\Theta^{\text{ext}} = \{(0.84,
0.79), (0.79, 0.84)\}$ and so CS only has to compute 2 PMFs as opposed to P and
AO which have to compute 81. We initialise with $\Theta^{\rgf{1}} = \{\hat{p}\}
= \{(0.75, 0.75)\}$. In iteration \rgf{1}, CS solves the MIP reformulation over
$\Theta^{\rgf{1}}$ and finds $y^{\rgf{1}} = 10$. It then evaluates the expected
costs under each $p \in \Theta^{\text{ext}}$ and finds the worst-case to be
given by $p^{\rgf{1}} = (0.84, 0.79)$. Hence, we have $\Theta^{\rgf{2}} =
\{(0.75, 0.75), (0.84, 0.79)\}$. In iteration 1, CS solves the model over
$\Theta^{\rgf{2}}$ and finds $y^{\rgf{2}} = 8$. It finds the worst-case cost to
be given by $p^{\rgf{2}} = (0.79, 0.84)$, and hence takes $\Theta^{\rgf{3}} =
\Theta^{\rgf{2}} \cup \{(0.79, 0.84)\}$. In iteration {\rgf{3}}, CS finds
$y^{\rgf{3}} = 9$ and $p^{\rgf{2}} = (0.84, 0.79)$. Since $(0.84, 0.79) \in
\Theta^{\rgf{3}}$, the algorithm ends and returns $y^{\text{CS}} = 9$ and
$p^{\text{CS}} = (0.84, 0.79)$ with an expected cost of $z^{\text{CS}} = 19.07$.
Hence, CS returned the optimal $y$ but slightly underestimated its worst-case
cost. This is an example of where CS will be suboptimal because $p^{\text{P}}
\notin \Theta^{\text{ext}}$. However, CS returned its solution in 0.17 seconds,
as opposed to P's 2.6 seconds. Note that CS terminated in 2 iterations because
$|\Theta^{\text{ext}}| \ = 2 = L$. 

To solve this model with AO, we construct the reduced set of intakes
$\tilde{\m{I}}$. In order to do so, we compute the PMFs, which takes 2 seconds.
Using $\beta = 0.001$, we find the new set of intakes to have $|\tilde{\m{I}}| \
= 150$, which is a 67\% cardinality reduction. Then, we solve the MIP model over
$\tilde{\m{I}}$ and find the solution $y^{\text{AO}} = 9$, $p^{\text{AO}} =
(0.82, 0.82)$, meaning that AO was both $y$-optimal and $p$-optimal in this
instance. However, it took 0.71 \rgf{seconds} in total, as opposed to CS's 0.17
seconds. We can also run this instance with CS\_opt. Doing so, CS\_opt's first
two iterations are the same as CS's. In its third iteration it finds
$y^{\rgf{3}} = 9$ and $p^{\rgf{3}} = (0.82, 0.82)$, whereas CS found
$p^{\rgf{3}} = (0.84, 0.79)$. Following this, in iteration $k=\rgf{4}$ it finds
$p^{\rgf{4}} =  (0.82, 0.82)$ and breaks since $p^{\rgf{4}} \in
\Theta^{\rgf{4}}$, returning $y^{\text{CS\_opt}} = 9$ and $p^{\text{CS\_opt}} =
(0.82, 0.82)$. This is the same solution as P gave. This took CS\_opt a total of
0.43 seconds. It finished in twice as many iterations as CS. 

%% file: 4_experimental_design.tex
\section{Design of Computational Experiments}\label{sec:exp_des}

This section details our experiments evaluating the performance of the
algorithms described in Section~\ref{sec:CS_AO} in comparison with the solution
from the parametric DRO model. These experiments will also allow us to compare
the solutions resulting from the \rgf{parametric model (P)} and the
\rgf{non-parametric model (NP)} and the times taken to reach optimality by each
model. In this section, we discuss how the parameters for the experiments will
be chosen to ensure that they are representative of typical real-life scenarios.
To discuss experimental design, we need to define which parameters of the model
will be varied and the values that they will take. The vector of inputs to the
model for a fixed set $\m{I}$ of intakes and $\m{P}$ of distributions is $S =
(c, D, a, L, K)$.

\subsection{Parameter Hierarchy}\label{sec:hierarchy}

It is helpful to consider a hierarchy of parameter choices, which is defined by:
\begin{enumerate}
    \item $(L, K)$ defines the difficulty of the problem in terms of the MIP
    itself.
    \item $c$ and $D$ define the set of solutions that are possible for a given
    model with fixed $L$ and $K$\rgf{. They} need to be constructed for each
    combination of \rgf{$L$ and $K$} to ensure \rgf{that} we have a varied range
    of instances when it comes to pulling forward opportunities. We create this
    variety by varying the number of days that have spare capacity and are hence
    able to receive additional jobs. The values of $c$ and $D$ used are
    discussed in Section~\ref{sec:c_D}.
    \item 
        \begin{enumerate}
            \item For the parametric model, $\m{I}$ and $\Theta$ define how the
            uncertainty is encoded in the model, depending on the planner's
            attitude to risk. If $|\m{I}|$ or $|\Theta|$ \rgf{is} large, solving
            to optimality will be very slow, and we would like to use a
            heuristic that is not significantly affected \rgf{by these sizes}.
            $|\m{I}|$ is defined by $i^{\max}$, and $|\Theta|$ is defined by two
            parameters. The initial discretisation of the interval $[0,1]$ in
            which each $p_\tau$ lies is defined by $n_{\text{probs}}$. The
            maximum distance from the nominal distribution that $p \in \Theta$
            can be is defined by the second parameter, $N$. This is the number
            of samples that we take from the distribution of $I$ in order to
            calculate $\hat{p}$. Larger $N$ results in smaller distances from
            $\hat{p}$ being allowed, and hence corresponds to a less risk-averse
            planner. For these experiments, we use $95\%$ confidence sets, i.e.\
            $\Theta_\alpha$ from~\eqref{eq:conf_set_fast} with $\alpha = 0.05$.
            From now on, we use $\Theta$ to represent $\Theta_{0.05}$.
            \item For the non-parametric model, we also use $95\%$ confidence
            sets. However, for this model we use the $\phi$-divergence based
            set, $\m{P}_\rho$ given in~\eqref{eq:NP_conf_set} with $\rho$
            defined by~\eqref{eq:rho} \rgf{and} $\alpha=0.05$. This set is only
            affected by $N$, which affects the maximum distance from $\hat{P}$
            that a distribution can lie under the non-parametric model.
        \end{enumerate}
    \item $a$ will be left as the ones vector for these experiments as it has
    not been seen to have an effect on solutions.
\end{enumerate}
We choose $L=5$ due to it being the number of days in a typical working week. We
take the maximum pulling forward window length to be  $K=2$\rgf{. This is
because} pulling forward is not enacted until the operational planning phase,
where the planning horizon is very short. These choices are partly motivated by
usual practices, and also partly by the \rgf{following} fact\rgf{. We} aim to
test our heuristics against optimal solutions, and for larger $L$ or $K$ the
model becomes very difficult to solve to optimality.  Note that the optimality
tolerance for CS/CS\_opt, $\varepsilon$, will be set \rgf{to} 0.01 and it will be
\rgf{run} for a maximum of $\rgf{k^{\max}} = 10$ iterations. Initial testing
suggested \rgf{that} these parameters are not so important, as CS and CS\_opt
\rgf{always} terminated due to a repeat parameter (i.e.\ $p^k \in \Theta^k$) after \rgf{less than 10} iterations.

\subsection{Capacity and Workstacks}\label{sec:c_D}

The factors affecting the potential solutions of a model the most are $c$ and
$D$, due to the fact that they define the rollover and pulling forward
opportunities. In this section, we detail the capacities $c$ and workstacks $D$
used in our experiments. These are constructed with the aim of ensuring that a
variety of combinations of pulling forward opportunities are represented by at
least one $(c, D)$ pair. We assume for this section that the previous parameters
in the hierarchy, i.e.\ $L$ and $K$\rgf{,} are given. We now define how $c$ and
$D$ define pulling forward opportunities mathematically. Firstly, we define the
set of pairs of days under consideration for pulling forward as:
\begin{equation}
    \m{F} = \l\{(\tau_1,\tau_2)\ \big| \ \tau_1 \in \{2,\dots,L\}, \tau_2\in
    \{\tau_1-K,\dots,\tau_1-1\}\r\}
\end{equation}
and the set of pairs such that the corresponding $y$ can feasibly be positive
given $c$ and $D$ as:
\begin{equation}
    \m{F}^+(c, D) = \l\{(\tau_1,\tau_2) \in \m{F}\ | \ c_{\tau_2} > D_{\tau_2},
    D_{\tau_1} > 0\r\}.
\end{equation} 
This is the set of all pairs of days $(\tau_1, \tau_2)$ such that $\tau_2$ is
within pulling forward range of $\tau_1$, $\tau_2$ has spare capacity and
$\tau_1$ has workstack jobs to be completed early. For our experiments, we
consider instances where $D_{\tau} > 0$ for all $\tau \in \{1,\dots,L\}$. This
is because for a short horizon of $L=5$ days, it is very unlikely that any day
will have a workstack of zero. Hence, we can control $|\m{F}^+(c, D)|$ by
controlling which days have spare capacity. For example, we can set $|\m{F}^+(c,
D)| \ = 3$ by setting $D_1 < c_1$ and $D_4 < c_4$ and then $D_\tau > c_\tau$ for
$\tau \in \{2, 3, 5\}$. This results in $\m{F}^+(c, D) = \{(2, 1), (3, 1), (5,
4)\}$. We do this similarly for other values of $|\m{F}^+(c, D)|$. The main
effect that $c$ and $D$ has on decision making is \rgf{that they define} the
constraints on $y$, meaning their only \rgf{important} quality is how much
pulling forward they do or do not allow. Using this set of values for $c$ and
$D$ we will be able to see how well our algorithms detect and make use of
opportunities for pulling forward.

\subsection{Uncertainty and Ambiguity Sets}

As a reminder, the term ``uncertainty set'' refers to $\m{I}$ and ``ambiguity
set'' refers to $\Theta$. We now detail the parameters used to construct these
sets in our instances.

\subsubsection{Uncertainty Set}

We assumed in Section~\ref{sec:problem_setting} that we would be given a set
$\m{I}$, either by expert knowledge or by a prediction model. We could then
extract $i^{\max}$ from this set. However, in these experiments, we do not have
access to real intake data or expert knowledge\rgf{. T}hus\rgf{,} it is more
convenient to define $i^{\max}$ and then use this to construct $\m{I}$. Since
there is a one-to-one mapping between the two, both methods achieve the same
result. We consider $i^{\max}$ satisfying:
\begin{equation}\label{eq:i_max_constraint}
    \sum_{\tau = 1}^L i^{\max}_\tau \le \sum_{\tau = 1}^L \max\{c_\tau - D_\tau,
    0\}.
\end{equation}
This is reasonable because if the total number of jobs arriving in the system
exceeds the RHS of~\eqref{eq:i_max_constraint} then some intake jobs will always
remain incomplete at the end of day $L$, regardless of our pulling forward
decision. Furthermore, we can \rgf{vary the} number of high-intake days, through
the \rgf{quantity} $n(i^{\max}) = \lvert\l\{\tau \in [L]: i^{\max}_\tau > c_\tau
- D_\tau\r\}\rvert$\rgf{. This} corresponds to the number of days with \rgf{the}
potential \rgf{for} spikes in demand. Depending on $c$ and $D$, $n(i^{\max})$
can range between 0 and $L-1$. However, for these experiments we consider
$n(i^{\max}) \in \l\{1, \l\lfloor{\frac{L}{2}}\r\rfloor, L - 1\r\}$ for
sufficient coverage of cases. The case of \rgf{$n(i^{\max}) = 1$} corresponds to
a one-day spike caused by an event such as a major weather event. The case of
\rgf{$n(i^{\max}) = \l\lfloor{\frac{L}{2}}\r\rfloor$} could correspond to an
extended spike lasting for multiple consecutive days, for example, a network
problem causing lots of service devices to break. The final case of
\rgf{$n(i^{\max}) = L - 1$ corresponds to} $L-1$ small spikes in intake, marking
a period of consistently high intake.

\subsubsection{Ambiguity Sets}

The choice of parametric ambiguity set depends on the choice of discretisation
of $[0,1]^L$ and also the way we in which we then reduce its size. The choice of
discretisation is defined by the parameter $n_{\text{probs}}$, and increasing
this value increases the size of the ambiguity set. For these experiments, we
consider $n_{\text{probs}} \in \{5, 10, 15\}$. In our preliminary testing we
found that any value larger than 15 can lead to intractability when solving the
parametric model to optimality.

Both ambiguity sets are also defined by the sampling parameter $N$. For the
purpose of testing our models, we consider $N \in \{10, 50, 100\}$. Clearly,
higher $N$ leads to better convergence to the true distribution of the
MLE/$\phi$-divergence, but it also leads to much smaller ambiguity sets and
typically less conservative decisions. Even for $N=50$, we obtained some
singleton ambiguity sets. Typically, \rgf{$N$} would be chosen by the planner
who is in control of the sampling process. However, the results of our testing
can be used to understand the tradeoff between the accuracy of the approximation
and the conservativeness of the resulting decisions. Hence, it may influence the
value of $N$ used by the planner. In these experiments, we will assume $\hat{i}
= (0.75i^{\max}_1,\dots, 0.75 i^{\max}_L)$. Hence, we will obtain $\hat{p} =
(0.75,\dots,0.75)$. In practice, $\hat{p}$ would be obtained from sampling the
true intake distribution. However, without access to true intake data, we set
the value somewhat arbitrarily, since it is only used for testing purposes. If
these models were used by a real planner, we would suggest that they
\rgf{calculate} their own MLE.

%% file: 5_results.tex
\section{Results}\label{sec:results}

We now detail the results of our experiments that \rgf{we used to} test the
algorithms on 279 problem instances \rgf{with} $L=5$ and $K=2$. We report the
results from all 5 algorithms in terms of times taken, pulling forward decisions
and worst-case distributions. Due to space considerations, we present some
additional results in the Appendices. We discuss the effects of workstacks on
solutions in Appendix B.1. We give a brief comparison of our results with those
from the RO version of the model in Appendix B.2. In addition, we present and
test a Benders decomposition algorithm for this problem in Appendix C. These
experiments were run in parallel on a computing cluster (STORM) which has 486
CPU cores. The solver used in all \rgf{instances} was the Gurobi Python package,
gurobipy~\citep{gurobi}. The version of gurobipy used was 9.0.1. The node used
on STORM was the Dantzig node, which runs the Linux Ubuntu 16.04.6 operating
system, Python version 2.7.12, and 48 AMD Opteron 638 CPUs.

\subsection{Summary of Instances and Their Sizes}

In Table \ref{tab:set_sizes}, we summarise the sizes of the sets $\m{I}$ and
$\Theta$, that form\rgf{ed} the basis for \rgf{the} constraints and variables in
the model. A summary of the sizes of $\m{I}$ is given in
Table~\ref{tab:i_range_size}. The table shows 7 of the 31 $i^{\max}$ values
considered and the size of the resulting set $\m{I}$. The other $i^{\max}$
values considered were permutations of the values shown in the table, and hence
led to $|\m{I}|$ values that are already listed in the table.
Table~\ref{tab:p_set_size} shows the values of $N$ and $n_{\text{probs}}$ used
and the average size of the resulting ambiguity sets. The sizes vary as the
construction of the set also depends on $i^{\max}$. The \rgf{instances} where
$|\Theta|\ = 1$ correspond to \rgf{instances} where $\frac{\chi^2_{L,
1-\alpha}}{N}$ was too small to allow any $p$ other than $\hat{p}$ to be in the
ambiguity set defined by~\eqref{eq:conf_set_fast}.
\begin{table}[htbp!] 
    \begin{subtable}[t]{.45\linewidth}
        \centering
        \input{new_images_tables/i_sizes_short.tex}
        \caption{Example $i^{\max}$ values and sizes of the associated
        uncertainty sets $\m{I}$ considered.}
        \label{tab:i_range_size}
    \end{subtable}
    \begin{subtable}[t]{.45\linewidth}
        \centering
        \input{new_images_tables/fast_Theta_sizes}     
        \caption{Parameters defining ambiguity sets and average size of
        corresponding sets.}   
        \label{tab:p_set_size}
    \end{subtable}
    \caption{Summary of input parameters and corresponding set sizes}
    \label{tab:set_sizes}
\end{table}

\vspace{-1em}
We can see here that our choices of $i^{\max}$ gave instances with as many
distinct intakes (and rollover vectors) as 20000, and as few as 392. The sizes
of the ambiguity sets varied between 1 and 8854, where the largest sets resulted
from the smallest $i^{\max}$ and $N$ values, and the largest $n_{\text{probs}}$
values. This is because the criteria for $p$ being included in $\Theta$ was
\rgf{$\sum_{\tau=1}^L N i^{\max}_\tau\frac{(\hat{p}_\tau -
p_\tau)^2}{\hat{p}_\tau (1-\hat{p}_\tau)} \le \chi^2_{L, 1-\alpha}.$} Clearly
the LHS is increasing in $N$ and $i^{\max}_\tau$. Hence, larger values lead to a
higher distance from the nominal distribution. Large $n_{\text{probs}}$ leads to
larger $\Theta$ because it results in a finer discretisation of $[0,1]^L$, and
hence more candidate $p$ values.

\subsection{Optimality of Algorithms and Times Taken}

Comparing results for DRO problems is not as simple as comparing final objective
values. Our optimal objective value can be written as $z^* = \min_y \max_p f(y,
p)$. Here\rgf{,} $f(y, p)$ is the total expected rollover cost, i.e.\
$\sum_{\tau=1}^L a_\tau\E_p(R_\tau\ |\ y)$. Suppose we have an instance where
$y^{\text{CS}} = y^{\text{P}}$ but $p^{\text{CS}} \neq p^{\text{P}}$. Then, if
CS gives a lower objective value than P, it may appear to have given a better
solution to the minimisation problem. However, this means that CS did not
successfully choose the worst-case $p$ for its chosen $y$. This leads to a lower
objective function value but a suboptimal solution \rgf{with respect to}\ $p$.
Similarly, we can say that CS is suboptimal if $p^{\text{P}} = p^{\text{CS}}$
but $y^{\text{CS}} \neq y^{\text{P}}$ and CS gave a higher objective value.
Hence, both a higher and a lower objective value can suggest suboptimality for a
DRO model. Given this, we summarise the results using 3 optimality criteria. An
algorithm $x \rgf{\ \in \{\text{CS}, \text{CS\_opt}, \text{AO}\}}$ is said to
be:
\begin{enumerate}
    \item \textbf{$y$-optimal} if $\max_{p\in \Theta} f(y^x, p) = z^*$.
    \item \textbf{$p$-optimal} for a given $y^x$ if $f(y^x, p^x) = \max_{p\in
    \Theta} f(y^x, p)$.
    \item \textbf{Optimal} if $f(y^x, p^x) = z^*$. Note that this is met is the
    algorithm is both $y$-optimal and $p$-optimal.
\end{enumerate}
We display the number of times each algorithm was optimal, $p$-optimal and
$y$-optimal in Table~\ref{tab:optimality}. 
\vspace{-1em}
\begin{table}[htbp!]
    \centering
    \input{new_images_tables/final2_optimality_counts}
    \caption{Summary of optimality of heuristics}
    \label{tab:optimality}
\end{table}

\vspace{-1em}
Table~\ref{tab:optimality} shows that CS was optimal in 92\% of instances, and
$y$-optimal in 97\%. As can be expected, CS\_opt was optimal in every instance.
AO was only optimal in 80\% of instances and $y$-optimal in 86\% of instances.
In fact, both CS and AO were optimal in selecting $p$ in more than 92\% of
instances. Unsurprisingly, AO performs the best in this regard. This is because
it solves the problem over the full set of parameters, unlike CS. However, CS
was still $p$-optimal in around 92\% of \rgf{instances}. Closeness to optimality
of the algorithms is discussed in Section~\ref{sec:results_summary}.

A summary of the computation times of each algorithm is given in
Table~\ref{tab:times_taken}. Firstly, the table shows average and maximum times
taken over all instances. CS \rgf{took} around 17 seconds \rgf{on average}. To
find the optimal solution, it took approximately 1 minute and 50 seconds on
average when using P, which is a large difference. CS\_opt found the optimal
solution in an average of 20 seconds, which is faster than P. This is only 3
seconds slower than CS on average. However, there are many instances with small
ambiguity sets. AO \rgf{took} similar times to CS\rgf{; it also took} around 17
seconds on average. NP solve\rgf{d} faster than P, but slower than CS, CS\_opt
and AO. The fact that NP \rgf{was} slower than CS\_opt suggests that the
parametric model can be solved to optimality faster than the non-parametric
model. AS reports the times taken to compute $\Theta$ for the parametric
algorithms. This \rgf{was} not included in the solution time for each algorithm,
as it is a pre-computation step. It is worth noting that the average of 6
seconds is significantly faster than extracting $\Theta$ from the non-parametric
confidence set, which can take hours. Please note that, while the differences
between the algorithms' times may seem small, these instances are small compared
to real planning instances. We would expect the time differences to be more
pronounced when the problems are large. Furthermore, CS\_opt requires
significantly more memory and computing power than CS. For \rgf{instances with
large ambiguity sets}, it stores thousands of distributions, each of which
comprises thousands of values. CS only stores around \rgf{$L$} distributions,
regardless of \rgf{the size of the ambiguity set}.
\begin{table}[htbp!]
    \centering
    \input{new_images_tables/final2_times}
    \caption{Summary of times taken}
    \label{tab:times_taken}
\end{table}

\vspace{-1em}
Since there \rgf{were} a large number of small instances that affect\rgf{ed} the
overall averages, Table~\ref{tab:times_taken} also shows average and maximum
times for instances with the largest ambiguity sets. This corresponds to the
largest 10\% of instances with respect to $\Theta$ or equivalently $|\Theta| \
\ge 1000$. From these two columns, we see that CS\_opt \rgf{took} more than 4
times longer than CS on average, when $\Theta$ \rgf{was} large. We also see that
CS\_opt took 34 seconds longer to solve its slowest instance than CS took for
its slowest instance. The largest time difference was 46 seconds, and this
occurred when $|\Theta|\ = 831$ and $|\m{I}| \ = 20000$. Th\rgf{is} time
difference \rgf{was} due to two main reasons. Firstly, CS never spent more than
0.5 seconds computing PMFs, whereas CS\_opt took up to 22 seconds. Hence, CS
significantly reduce\rgf{d} the amount of computation required. Secondly, CS
typically completed in many fewer iterations than CS\_opt. This is because its
use of $\Theta^{\text{ext}}$ meant it identified a repeat parameter in fewer
iterations. Based on the optimality counts and time taken, CS is the strongest
heuristic. It select\rgf{ed} the optimal $y$ is 97\% of instances, and \rgf{did}
so in less time than CS\_opt. CS\_opt can be used when $\Theta$ is small, but it
\rgf{will begin} to solve slowly in comparison with CS when $\Theta$ is large.

\subsection{Performance of Algorithms in Detail}\label{sec:results_summary}

To illustrate further how well the algorithms performed, we define the following
two metrics. Note that a positive value for \textit{either} of these metrics
suggests suboptimality.
\begin{enumerate}
    \item \textbf{Quality of $p$ choices.} For a solution $y^x$ that was
    selected by an algorithm $x$, where $x \in \{\text{CS}, \text{CS\_opt},
    \text{AO}\}$, we calculate the worst-case expected cost over all
    distributions in $\Theta$ using brute force. We can then compare this cost
    with the expected cost obtained by the algorithm, i.e.\ from $p^x$, the $p$
    that the algorithm selected. This allows us to establish how close to
    worst-case the choices in $p$ were. We refer to this difference as the
    $p$-gap, and it is defined as $g_{p}(y^x, p^x) =  \max_{p \in \Theta} f(y^x,
    p) - f(y^x, p^x).$
    \item \textbf{Quality of $y$ choices.} For a given solution $y^x$ from
    algorithm $x$, we compute the worst-case expected cost using brute force, as
    we did when finding $g_{p}(y^x, p^x)$. We can then compare this worst-case
    cost with that of the optimal $y$, to assess how close $y^x$ is to optimal.
    This is referred to as the $y$-gap, and is defined as $g_{y}(y^x) =
    \max_{p\in \Theta}f(y^x, p) - z^*.$
\end{enumerate}
In Table~\ref{tab:gap_summaries}, we summarise the average $p$-gaps and $y$-gaps
of the three heuristics, along with the average absolute percentage gaps (APGs).
The $p$-APG \rgf{was} obtained by taking the $p$-gap as an absolute percentage
of the worst-case expected cost for the chosen solution $y^x$. The $y$-APG
\rgf{was} obtained by taking the $y$-gap as an absolute percentage of the
optimal objective value.
\begin{table}[htbp!]
    \centering
    \input{new_images_tables/final2_gaps}
    \caption{Summary of gaps and APGs of the heuristics}
    \label{tab:gap_summaries}
\end{table}

\vspace{-1em}
This suggests that all algorithms perform very well at choosing the worst-case
$p$ for a fixed $y^x$, \rgf{since} all \rgf{had} an average $p$-APG of less than
0.09\%. AO perform\rgf{ed} the best at selecting $p$, which supports the
observation made from the optimality counts. CS and AO are very good at
selecting the optimal $y$, since they both use a solver to do so. Of CS and AO,
CS perform\rgf{ed} the best in this regard, with an average $y$-APG of 0.01\%.
The $y$ solution CS chose had, on average, a worst-case expected cost that was
0.0058 away from the optimal objective value. AO also perform\rgf{ed} well in
selecting $y$, but its average $y$-APG \rgf{was} a factor of 20 larger than that
of CS. Due to its optimality in every instance, CS\_opt ha\rgf{d} average gaps
and APGs of 0.

\begin{figure}
    \centering
    \begin{subfigure}[t]{0.45\textwidth}
        \includegraphics[width=0.9\textwidth]{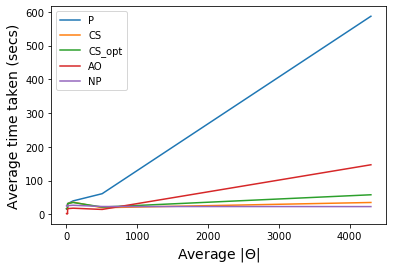}
        \caption{Average times taken by $|\Theta|$}
        \label{fig:times_by_theta}    
    \end{subfigure}
    \begin{subfigure}[t]{0.45\textwidth}
        \includegraphics[width=0.9\textwidth]{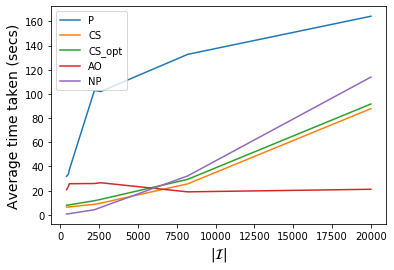}
        \caption{Average times taken by $|\m{I}|$}
        \label{fig:times_by_I}    
    \end{subfigure}
    \caption{Average times taken by sizes of sets}
    \label{fig:times_by_set}
\end{figure}

We also study the results broken down by the size of the set of distributions.
In order to reduce the size of the table, we present results averaged over the
categories for $|\Theta|$ given in Table~\ref{tab:p_set_size}. We present these
results in Table~11, which is in Appendix~D.1 due to space considerations. In
summary, the table suggests that CS \rgf{did} not return suboptimal $p$s for its
chosen $y$ until the set reache\rgf{d} the average size of 93. CS \rgf{was}
consistent in its selection of the optimal $y$ across all values of $|\Theta|$.
CS's $y$-APE stay\rgf{ed} very close to 0 in all \rgf{instances}. AO ha\rgf{d}
larger $p$-gaps for the larger values of $|\Theta|$. Interestingly, AO's
performance in selecting $y$ improve\rgf{d} as  $|\Theta|$ \rgf{grew} larger.
CS\_opt ha\rgf{d} zero gaps and APGs for all values of $|\Theta|$, but its times
taken d\rgf{id} not scale as well as CS's and AO's with large $|\Theta|$. For
small ambiguity sets, CS\_opt \rgf{took} similar times to CS, but it \rgf{took}
twice as long for the largest ambiguity sets (average size of 4301). We also
plot the average times by $|\Theta|$ in Figure~\ref{fig:times_by_theta}.
\rgf{This plot suggests that the} algorithms that use Gurobi on the full set of
distributions, i.e.\ P and AO, do not scale well with $|\Theta|$ in terms of
time. CS, CS\_opt and NP all scale much better with $|\Theta|$ than AO and P.
For CS and CS\_opt, this is because they only ever solve an MIP reformulation
over a small subset of $\Theta$. For NP, this is because increasing the size of
the ambiguity set for the non-parametric model does not result in a more complex
model, it only increases the value $\rho$. This plot supports our
\rgf{conclusion} that CS\_opt solves in similar times to CS when $\Theta$ is
small, but begins to take noticeably longer for large $\Theta$. 

Finally, we can look at the performance of \rgf{the} algorithms by the size of
the set of intakes $\m{I}$. These results are shown in Table~12 in Appendix D.2.
The $p$-APGs for the heuristics \rgf{were} not significantly affected by
$|\m{I}|$, apart from a drop in performance for CS when $|\m{I}| \ = 8192$. This
\rgf{was} likely due to other model parameters, since there is no reason for
$|\m{I}|$ to affect the $p$-APG. AO also beg\rgf{an} to lose $y$-performance
when $|\m{I}|\ = 8192$. This is an intuitive result, because as this set gets
larger AO will remove more and more intakes. This reduces the accuracy of its
approximation of the objective function. CS does not remove intakes, which
explains why its performance \rgf{was} consistent. In fact, CS's $y$-APG
\rgf{was} lower than AO's when $|\m{I}| \ = 20000$. Again, CS\_opt ha\rgf{d} all
zero gaps and APGs. The difference between CS and CS\_opt in terms of times
taken is less noticeable here. CS\_opt consistently \rgf{took} 3-5 seconds
longer than CS for all values of $|\m{I}|$. This indicates that $|\Theta|$
\rgf{was} the main factor causing CS\_opt to solve slowly. We also plot the
average times by $|\m{I}|$ in Figure~\ref{fig:times_by_I}. This plot suggests
that P does not scale well with $|\m{I}|$, and that AO scales very well with
$|\m{I}|$. CS, CS\_opt and NP scale better than P, but not nearly as well as AO,
due to the fact that they do not apply dimension reduction to $\m{I}$.

\subsection{CS's Suboptimal Distributions} \label{sec:CS_subopt}

In this section, we compare the solutions and distributions from CS with those
from P. Since CS is only limited by its performance in selecting $p$, we study
CS's worst-case $p$s in order to find ways to improve its performance. We do not
study CS's performance with respect to $y$, since if $\Theta^{\text{ext}}$
contains $p^{\text{P}}$ then CS will return the same $y$ as P, as evidenced by
CS\_opt. Hence, improving $\Theta^{\text{ext}}$ is sufficient to improve CS with
respect to $y$ and $p$. We do not analyse AO's solutions, since improving its
performance can only come from tuning $\beta$.

As shown in Table~\ref{tab:optimality}, CS chose the optimal $p$ for its
selected $y$ in 92\% of \rgf{instances}, leaving 22 instances where it did not.
This indicates that our set $\Theta^{\text{ext}}$ did not in fact contain the
worst-case $p$ in those 22 instances. To compare CS with P, we study only
instances where CS selected the same $y$ as P, which occurred in 15 of these 22
instances. Firstly, for these 15 instances, we can confirm that $p^{\text{P}}$
was not contained in the set $\Theta^{\text{ext}}$ used by CS. This either
occurred because no probability was at its maximum, or because the sum of the
probability vector was not maximised. We find that one value of $p^{\text{P}}$
was maximised in 13 out of 15 instances. However, in every one of these 13
instances, the sum over the vector was not maximised. This indicates that the
main reason why CS did not return the worst-case $p$ in every \rgf{instance} was
because the worst-case does not need to satisfy this condition. In general, we
find that CS both allocated a higher maximum success probability and more
success probability in total than P.

In order to see why the worst-case $p$ does not need to satisfy the
sum-maximisation criterion, we study some examples more closely. For example, in
one instance we \rgf{had} $p^{\text{P}} = (0.933, 0.867, 0.867, 0.867, 0.733)$
and $p^{\text{CS}} = (0.933, 0.933, 0.867, 0.8, 0.867)$. We see that P and CS
both gave maximal probability to day 1. However, P reduced days 2 and 5's
probabilities in order to allocate more to day 4. The resulting rollover vectors
\rgf{were} $(0.87, 10.6, 27.34, 24.54, 41.01)$ for P and $(0.87, 10.74, 27.47,
24.27, 41.0)$ for CS. In this instance, allocating higher probability to day 4
resulted in higher day-4 and also day-5 rollover, and more rollover in total,
despite the fact that the total probability was not maximised. Another way that
CS can be suboptimal is choosing the wrong day to set to its maximum. For
example, in one of the 15 instances P gave $p^{\text{P}} = (0.933, 0.867, 0.8,
0.8, 0.8)$ and CS gave $p^{\text{CS}} = (0.8, 0.867, 0.8, 0.867, 0.867)$. Here,
CS set $p_2$ to its maximum, while P set $p_1$ at its maximum. For this
instance, the closest values of $p$ to $p^{\text{P}}$ that were in
$\Theta^{\text{ext}}$ were $(0.933, 0.8, 0.8, 0.867, 0.933)$ and $(0.933, 0.8,
0.8, 0.933, 0.867)$. These two solutions give less expected cost than
$p^{\text{CS}}$, and so CS did not allocate maximal probability to day 1.
Clearly the maximal cost came from allocating high probability to day 2 as well
as day 1, but no such probability vectors were contained in
$\Theta^{\text{ext}}$. Finally, in two instances no value of $p^{\text{P}}$ was
at its maximum. One example of this occurred when $p^{\text{P}} = (0.867, 0.867,
0.867, 0.733, 0.733)$ and $p^{\text{CS}} = (0.8, 0.933, 0.8, 0.8, 0.8)$. CS has
allocated day 2 its maximum probability. However, the worst-case parameter
spread the success probability more evenly over the first 3
days.% The assumption that one $p_\tau$ is maximised is usually satisfied, and hence this is not the main reason for CS's suboptimality.

These observations explain why CS did not always return the true worst-case $p$.
Clearly, the issue lies in the construction of  $\Theta^{\text{ext}}$. In
particular, the assumption that the sum over the success probability vector
should be maximised is not always required. In fact, sometimes it is worse to
reduce the sum in order to give high priority days a higher success probability.
In order to assess whether or not this is the case, CS would need to compare the
maximum intakes for each day in order to see where the most rollover could be
caused. 

\subsection{Parametric vs.\ Non-parametric Decisions and
Distributions}\label{sec:P_vs_NP}

In this section, we compare NP's solutions and distributions with those from P.
This will allow us to assess the benefits and costs of including the parametric
information in the model. As we have seen, incorporating this information
creates a model that is larger and computationally more difficult to solve.
However, it retains the information on the family of distributions that $P^0$
lies in and ensures that the worst-case distribution from the model is also in
this family. This is something that is not guaranteed by the non-parametric
model.

\subsubsection{Pulling Forward Decisions and Objective Values}\label{sec:PvNP_y}

We first study the differences in pulling forward decisions between the two
models along with their worst-case objective values. We find that the two models
gave the same pulling forward decision in 199 of the 279 instances solved. This
can be stated as NP being $y$-optimal with respect to the parametric model in
71\% of instances. In every one of these instances, it was only optimal to pull
forward between either days 2 and 1 or not at all. The worst-case expected cost
from NP was $1.21$ higher than that from P in these instances, on average. This
suggests that the worst-case distribution from NP for a fixed $y$ is typically
worse than that from P.

\begin{figure}[htbp!]
    \centering
    \begin{subfigure}[t]{0.45\textwidth}
        \includegraphics[width=0.9\textwidth]{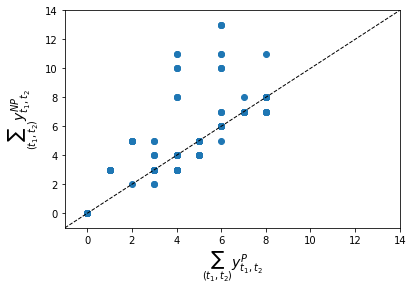}
        \caption{\rgf{Amount pulled forward under P vs.\ NP}}
        \label{fig:PvsNP_y_sum}        
    \end{subfigure}
    \begin{subfigure}[t]{0.45\textwidth}
        \centering
        \includegraphics[width=0.9\textwidth]{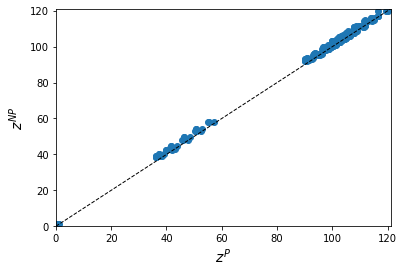}
        \caption{\rgf{Worst-case expected costs under P vs.\ NP}}
        \label{fig:PvsNP_obj}
    \end{subfigure}
    \caption{Scatter plots comparing P and NP's pulling forward decisions}
    \label{fig:PvsNP_y}
\end{figure}

\vspace{-1em}
Figure~\ref{fig:PvsNP_y_sum} shows a scatter plot of the total amount pulled
forward under each model in each of the 279 instances.
Figure~\ref{fig:PvsNP_obj} shows the corresponding worst-case expected costs.
The dashed line corresponds to instances where both models pulled forward the
same amount or had the same worst-case cost. The points in
Figure~\ref{fig:PvsNP_y_sum} where the decisions were different suggests that
there is no definitive answer to which \rgf{model's decision} is more
conservative. In 42 instances NP pulled forward more, and in 38 instances it
pulled forward less. However, when NP pulled forward more than P, it pulled
forward up to 7 jobs more. When P pulled forward more, it only pulled forward 1
\rgf{job} more. On average over the instances where the two solutions were
different, NP pulled forward 1.24 more jobs. The overall average difference was
0.32. This suggests that NP is generally slightly less conservative than P.
However, as shown in Figure~\ref{fig:PvsNP_obj}, rarely did NP attain a lower
worst-case expected cost than P. The overall average difference between P and
NP's worst-case expected costs was $-1.21$. This suggests that NP's worst-case
distribution typically suggests that there will be 1.21 more jobs being expected
to roll over in the worst case. This is surprising since NP typically pulled
forward more. Hence, this result indicates that NP's less conservative nature
\rgf{led} to more expected rollover in the majority of these instances.

Since NP results from relaxing the requirement that the worst-case distribution
is binomial, we can view NP as a heuristic for solving the parametric model.
Hence, it may be beneficial to study the expected cost resulting from
$y^{\text{NP}}$ under the binomial worst-case distribution that would be given
by P, instead of the distribution given by NP. Therefore, for each value of
$y^{\text{NP}}$, we compute the worst-case binomial distribution given by a $p
\in \Theta$, and the associated expected cost. This allows us to compute the
objective value that $y^{\text{NP}}$ would attain under the parametric model.
Hence, it allows us to assess the quality of $y^{\text{NP}}$ in comparison with
$y^{\text{P}}$, as we did for our heuristics. We can also study the difference
between $y^{\text{NP}}$'s worst-case cost under P and NP, via the $p$-gap. This
allows us to assess how the two objective functions differ for the same $y$.
\rgf{As a reminder, for an algorithm $x$ the $p$-gap is defined as $g_p(y^x,
p^x) = \max_{p \in \Theta} f(y^x, p) - f(y^x, p^x)$, and the $y$-gap is given by
$g_{y}(y^x) = \max_{p\in \Theta}f(y^x, p) - z^*$}.
\begin{table}[htbp!]
    \centering
    \input{new_images_tables/final2_NP_gaps}
    \caption{Summary of NP and CS's gaps}
    \label{tab:NP_gaps}
\end{table}

\vspace{-1em}
The gaps for NP are summarised in Table~\ref{tab:NP_gaps}, along with those from
CS for comparison. The $p$-gaps \rgf{show} that the worst-case cost for
$y^{\text{NP}}$ from NP \rgf{was} 1.18 higher than that from P, on average. This
indicates that the NP model typically overestimate\rgf{d} the worst-case cost
associated with $y^{\text{NP}}$. This is consistent with our previous
observation that NP's worst-case objective values \rgf{were} higher for a fixed
$y$. In fact, NP overestimated the worst-case cost of $y^{\text{NP}}$ in 248 of
279 instances (89\%). The most that NP overestimated this cost by was 3.4. These
values may seem small, but relative to the true worst-case cost they can be
quite large. The largest $p$-APG was 165\%, indicating that the worst-case cost
from NP was 2.65 times that from P. These results indicate that NP will
typically give an objective that makes a decision look worse than it would be in
reality.  The $y$-APGs suggest that the $y$ decisions from NP perform\rgf{ed}
similarly to that of P, under P's objective. However, they \rgf{did} result in a
slight cost increase on average. Based on the results here, we believe that CS
is the strongest performing algorithm. CS \rgf{ran} in less time than NP and
\rgf{gave} solutions closer to those from P. In fact, we can say that the NP
solutions \rgf{had} gaps that \rgf{were} 4 times higher than CS's on average.
Both average gaps \rgf{were} small, but CS was optimal in 92\% of instances, as
opposed to 71\% for NP. Furthermore, if one were to use the NP model, then they
would likely overestimate the rollover cost from their decision by approximately
13\%, whereas CS would underestimate this cost by approximately 0.084\%. 

\subsubsection{Worst-case Distributions}\label{sec:PvNP_p}

In order to explain the differences in decisions and costs, we now study the
worst-case distributions from P and NP. There are a number of ways in which
these distributions can be different. The most obvious one is that P's
worst-case distribution is always binomial, whereas NP's is not. As well as
this, the two approximations of the 95\% confidence set for $P^0$ can be
different, allowing different distances from $\hat{P}$.  In fact, typically the
confidence sets for P were larger. This indicates \rgf{that the parametric sets
had better coverage}. We first study the maximum distances from $\hat{P}$
allowed by each ambiguity set and the distances attained by the parametric and
non-parametric worst-case distributions, as measured by $d_\phi$. We find that
the maximum distance allowed by P \rgf{could} be almost twice that allowed by
NP. The maximum distance allowed by NP was 1.22, whereas this value was 2.32 for
P. This suggests that the parametric ambiguity set can be significantly larger
than the non-parametric set. We also find that NP's worst-case distribution
always achieve\rgf{d} the maximum distance from $\hat{P}$. Interestingly, the
same does not apply for P. The maximum distance that $P^{\text{P}}$ had from
$\hat{P}$ was 2.00, showing that the worst-case binomial distribution was not
always as far from $\hat{P}$ as it was allowed to be. In fact, there were 106
instances where P did not reach its maximum distance. As a result, even though
the parametric ambiguity sets allowed $P^{\text{P}}$ to be further from
$\hat{P}$, we still find that $P^{\text{NP}}$ was further from $\hat{P}$ on
average. The fact that NP's solution \rgf{was} always on the boundary may
indicate that the true worst-case distribution \rgf{was} further from $\hat{P}$
than \rgf{was} allowed by NP's ambiguity set. 

In order to compare the worst-case distributions directly, we compute a number
of summary statistics for each distribution and present their average values in
Table~\ref{tab:dists_summary}. \rgf{This table also shows the percentage
difference between the summary values for the two distributions, which is
calculated as $100\times\frac{\text{NP} - \text{P}}{\text{P}}$}. The first two
results we show are the average distances from $\hat{P}$ as measured by $d_\phi$
and by the Kullback-Leibler Divergence (KLD). The KLD value, $KLD(P^x,
\hat{P})$, can be loosely interpreted as the amount of \textit{surprise} that
would result in simulating from $P^x$ if the true distribution were $\hat{P}$.
These two rows indicate that NP was further from $\hat{P}$, on average, than P
with respect to both distance measures. The values of $d_\phi$ are quite close,
but proportionally the difference in $KLD$ values is much larger. In fact, NP
\rgf{had} 52\% more surprise than P, on average. This is likely due to the fact
that $P^{\text{NP}}$ is not binomial, unlike $P^{\text{P}}$. Entropy also
measures surprise, but with respect to the values given by the distribution. We
see that both distributions \rgf{had} a similar total entropy, but P \rgf{had}
slightly more. 

\begin{table}[htbp!]
    \centering
    \input{new_images_tables/dists_summary}
    \caption{Summary statistics comparing $P^{\text{P}}$ with $P^{\text{NP}}$}
    \label{tab:dists_summary}
\end{table}

\vspace{-1em}
We also present summaries of the total mean, variance and skewness of each
distribution. We see that NP \rgf{had} a higher total expected intake than P on
average, but less variance. This can be expected since NP can control the mean
and variance separately. P, on the other hand, fixes the variance by fixing the
mean. P can therefore have a smaller variance than P, even when the two means
are the same. However, P \rgf{was} typically less negatively skewed than NP.
These results may explain why NP's worst-case costs \rgf{were} higher. If NP is
more negatively skewed with a higher mean and lower variance, then this suggests
that more mass is allocated to the higher intakes and less to the lower ones.
Hence, expected costs will necessarily be higher.  

Finally, we look at the number of intakes that were \textit{popped} and
\textit{suppressed} by each worst-case distribution. A distribution $P^x$
popping an intake $i$ is defined as $\P(I = i \ | \ P^x) > 0$ when $\P(I = i\ |
\ \hat{P}) = 0$. The distribution $P^x$ suppressing $i$ is defined as $\P(I = i
\ | \ P^x) = 0$ when $\P(I = i\ | \ \hat{P}) > 0$. Since $\hat{P}$ is a binomial
distribution, technically we will never have popping as $\P(I = i\ | \ \hat{P})
> 0 \fa i \in \m{I}$. We will also never have suppressing under P, for the same
reason. In addition, by~\cite{Love15}, the modified $\chi^2$ divergence cannot
pop scenarios. Hence, we consider the distributions when rounded to 6 d.p.\
instead. The table shows that NP popped 55\% more intakes than P on average.
Both popped only a few intakes, which is consistent with our observation that
neither method can technically pop scenarios. This is not the main cause for the
difference in the distributions, however. The main difference is from
suppressing. We see that NP suppresse\rgf{d} 163\% more intakes than P on
average. This indicates that NP's worst-case distribution set a large number of
$\hat{P}$'s positive values to zero. P is much more restricted in this sense,
due to the fact that $P^{\text{P}}$ is also binomial. This means that P cannot
set any values to be exactly zero. This difference may also explain the
increased values of KLD given by NP; some intakes that would be generated by
$\hat{P}$ would not be generated by $P^{\text{NP}}$.

%% file: new_images_tables/i_sizes_short.tex
\begin{tabular}{lr}
\toprule
      $i^{\max}$ &  $|\mathcal{I}|$ \\
\midrule
 (1, 6, 6, 1, 1) &              392 \\
 (1, 3, 3, 3, 3) &              512 \\
 (2, 2, 2, 6, 2) &              567 \\
 (2, 2, 8, 8, 2) &             2187 \\
 (5, 5, 1, 5, 5) &             2592 \\
 (1, 7, 7, 7, 7) &             8192 \\
 (9, 9, 1, 9, 9) &            20000 \\
\bottomrule
\end{tabular}

%% file: new_images_tables/fast_Theta_sizes.tex
\begin{tabular}{rrr}
\toprule
 $N$ &  $n_{probs}$ &  Average $|\Theta|$ \\
\midrule
 100 &            5 &               1.000 \\
 100 &           10 &               1.000 \\
 100 &           15 &              16.871 \\
  50 &            5 &               1.419 \\
  50 &           10 &              14.419 \\
  50 &           15 &              93.129 \\
  10 &            5 &              14.742 \\
  10 &           10 &             504.226 \\
  10 &           15 &            4301.645 \\
\bottomrule
\end{tabular}

%% file: new_images_tables/final2_optimality_counts.tex
\begin{tabular}{llll}
\toprule
{} & No. (\%) Optimal Sol & No. (\%) $p$-Optimal Sol & No. (\%) $y$-Optimal Sol \\
\midrule
CS     &        257 (92.11\%) &            257 (92.11\%) &            272 (97.49\%) \\
CS\_opt &        279 (100.0\%) &            279 (100.0\%) &            279 (100.0\%) \\
AO     &        223 (79.93\%) &            263 (94.27\%) &            239 (85.66\%) \\
\bottomrule
\end{tabular}

%% file: new_images_tables/final2_times.tex
\begin{tabular}{lllll}
\toprule
{} & Avg. t.t. (Overall) & Max t.t. (Overall) & Avg. t.t. (Large) & Max t.t. (Large) \\
\midrule
P      &          0:01:22.85 &          0:19:23.5 &        0:07:57.99 &        0:19:23.5 \\
CS     &          0:00:17.48 &         0:01:50.37 &         0:00:06.1 &       0:00:33.82 \\
CS\_opt &          0:00:20.17 &         0:02:12.09 &        0:00:24.74 &       0:01:07.84 \\
AO     &          0:00:17.29 &         0:03:52.97 &        0:02:08.08 &       0:03:52.97 \\
NP     &          0:00:25.35 &         0:03:26.88 &        0:00:07.25 &        0:00:44.4 \\
AS     &          0:00:05.95 &         0:00:19.38 &        0:00:14.18 &       0:00:16.59 \\
\bottomrule
\end{tabular}

%% file: new_images_tables/final2_gaps.tex
\begin{tabular}{lrlrl}
\toprule
{} &  Avg. $p$-gap & Avg. $p$-APG &  Avg. $y$-gap & Avg. $y$-APG \\
\midrule
CS     &        0.0561 &      0.084\% &        0.0058 &     0.0101\% \\
CS\_opt &        0.0000 &        0.0\% &        0.0000 &        0.0\% \\
AO     &        0.0064 &     0.0065\% &        0.0233 &     0.1369\% \\
\bottomrule
\end{tabular}

%% file: new_images_tables/final2_NP_gaps.tex
\begin{tabular}{lrlrll}
\toprule
{} &  Avg. $p$-gap & Avg. $p$-APG &  Avg. $y$-gap & Avg. $y$-APG & $y$-opt. \% \\
\midrule
NP &       -1.1764 &    13.0095\% &        0.0234 &     0.0429\% &      87.1\% \\
CS &        0.0561 &      0.084\% &        0.0058 &     0.0101\% &      97.1\% \\
\bottomrule
\end{tabular}

%% file: new_images_tables/dists_summary.tex
\begin{tabular}{lrrl}
\toprule
{} &      P &       NP &     \% Gap \\
\midrule
$d_{\phi}(P^x, \hat{P})$ &    0.435 &    0.480 &   10.345\% \\
$KLD(P^x, \hat{P})$      &    0.167 &    0.254 &   52.096\% \\
Entropy              &    5.379 &    5.227 &   -2.826\% \\
Total EV             &   16.701 &   17.048 &    2.078\% \\
Total Variance       &    3.670 &    3.590 &    -2.18\% \\
Total Skewness       &   -4.274 &   -4.431 &    3.673\% \\
No.\ Suppressed      &  215.556 &  568.178 &  163.587\% \\
\bottomrule
\end{tabular}

%% file: 6_conclusions.tex
\section{Conclusions and Further Research}

In this paper, we presented parametric and non-parametric DRO models for a
workforce planning problem under a mixture of known and uncertain demand. We
developed heuristics to solve the parametric model, due to its poor scalability.
The general conclusions that we can make from our results are as follows. The
full model can be slow to solve to optimality using the MIP reformulation, i.e.\
using P. CS\_opt solves this model to optimality in a short time on average, but
begins to solve slowly when the ambiguity set is large. Our heuristics, AO and
CS, employ dimension reduction to the sets of intakes and distributions
respectively in order to solve the problem in significantly less time than P.
The main conclusion we make about these algorithms is that CS performs very
well, and takes a fraction of the time that P takes. However, we found that CS
can fail to select the worst-case success probability for its chosen pulling
forward decision due to its assumption that the total success probability should
be maximised. We compared the parametric and non-parametric solutions, and made
a number of conclusions. Namely, NP typically pulls forward more than P but it
overestimates the worst-case cost of a decision. Our results also suggst that
the NP distributions have higher means, lower variance and more negative
skewness. They also suppress many more intakes than P's distributions.

The main contribution that we have made to the existing DRO literature is the
new modelling framework of parametric DRO. \rgf{In real-world planning problems,
incorporating distributional ambiguity often results in intractable models.
Instead, data-driven estimates of the demand distribution or its parameters are
commonly used. However, this can lead to poor solutions when the estimates are
poor. Our methodology provides a way that parameter estimates can be utilised in while also hedging against cases when they are inaccurate. It allows
planners to build confidence sets around their estimates, that can be adjusted
to fit their level of risk aversion. For example, if the planner does not have
trust in their estimates then they can choose a large confidence level in order
to generate a larger and hence more risk-averse ambiguity set. In addition, our
use of parametric distributions instead of non-parametric ones means that the
worst-case distribution from our model is more explainable, since it can be
summarised by a small number of parameters. This distribution is also less
extreme, and less surprising given the estimated distribution. Furthermore,
using parametric distributions has allowed us to create fast algorithms for
solving the planning problem. This means that planners can incorporate
additional uncertainty without having to wait long periods of time for
solutions.}
%This framework allows us to enforce that the worst-case distribution from the
%DRO model belongs to the same family of parametric distributions as the true
%distribution. For an example problem, we also contribute a selection of
%algorithms resulting from incorporating information about the family in which
%the  distribution of the uncertain parameters lies. We have shown that this
%additional information can be very beneficial in constructing algorithms for
%the DRO model, since CS exploits the behaviour of the binomial distribution in
%order to produce solutions. This has not only allowed these algorithms to
%perform very well, but also to quickly produce solutions to a large and complex
%model. It has been noted in the literature that CS algorithms suffer runtime
%issues due to the complexity of the distribution separation problem, and our
%CS\_opt algorithm supports this. However, we have presented a heuristic CS
%algorithm that does not suffer from this problem.

There are a number of natural extensions to our work which would be of further
interest from a practical viewpoint. Firstly, we have considered a simplified
problem in which each job requires one unit of capacity to complete. This is not
typically the case in real life workforce planning. Adding more varied
completion times would be a clear next step in improving this model. Secondly,
the model considers the case where there is only one skill, and is equivalent to
assuming all workers can complete any job. In some scenarios this is not the
case, and the model could account for this by considering separate demand values
and decision variables for each skill. Thirdly, we have treated the capacity as
fixed and aimed to optimise its use. In some cases, if not all, however,
capacity can be manipulated in the tactical planning phase. For example, one can
order extra units of existing resources (overtime) or hire outside resources for
a cost (contractors). These ways to manipulate capacity (planning levers) will
form the basis of some of our future research. Finally, we have assumed in this
paper that the intakes are independent. Extending our model to account for
correlated intakes is a promising area for future work.

%% file: 7_appendices.tex
\begin{appendices}

    \setcounter{equation}{37}
    \setcounter{table}{6}
        
    \section{Derivation of CQP Reformulation of Non-parametric
    Model}\label{sec:NP_reformulation}
    
    \subsection{General Reformulation}
        
    For a $\phi$-divergence ambiguity set with nominal distribution $Q$, we can
    write the inner problem of the DRO model (1)-(8) as:
    \begin{align}
        \max_{P} & \sum_{\tau=1}^L a_\tau \E_P(R_\tau) \\
        \text{s.t. } & P_j \ge 0 \fa j = 1,\dots,n\\
        & \sum_{j=1}^n P_j = 1 \\
        & d_{\phi}(P, Q) \le \rgf{\rho}.
    \end{align}
    The Lagrangian of this model is given by:
    \begin{equation}
        L(P, \lambda, \nu) =  \sum_{\tau=1}^L \sum_{j=1}^n P_j a_\tau R^{i^j}_\tau +
        \lambda\l(\rgf{\rho} - d_{\phi}(P, Q)\r) + \nu\l(1 - \sum_{j=1}^n P_j \r) .
    \end{equation}
    The objective function of the dual problem is therefore:
    \begin{equation}
        g(\lambda, \nu) = \max_{P \ge 0} L(P, \lambda, \nu).
    \end{equation}
    Since $\rgf{\rho} > 0$ and $d_{\phi}(Q, Q) = 0 < \rgf{\rho}$ where $Q$ is a
    feasible choice of distribution, Slater's condition holds. Since the primal
    is concave, we have strong duality. We can hence write the objective of the
    dual of the inner problem as:
    \begin{align}
        \min_{\lambda \ge 0,  \nu} g(\lambda, \nu) &= \min_{\lambda \ge 0,  \nu}
        \max_{P \ge 0} \l\{\sum_{\tau=1}^L \sum_{j=1}^n P_j a_\tau R^{i^j}_\tau +
        \lambda\l(\rgf{\rho} - d_{\phi}(P, Q)\r) + \nu\l(1 - \sum_{j=1}^n P_j \r)
        \r\}\\
        &= \min_{\lambda \ge 0,  \nu}\l\{ \lambda \rgf{\rho} + \nu + \max_{P \ge
        0}\l(\sum_{j=1}^n P_j \sum_{\tau=1}^L a_\tau R^{i^j}_\tau - \lambda
        d_{\phi}(P, Q) - \nu \sum_{j=1}^n P_j\r) \r\}\\
        &= \min_{\lambda \ge 0,  \nu}\l\{ \lambda \rgf{\rho} + \nu + \max_{P \ge
        0}\l(\sum_{j=1}^n P_j \sum_{\tau=1}^L a_\tau R^{i^j}_\tau - \lambda
        \sum_{j=1}^n Q_j \phi\l(\frac{P_j}{Q_j}\r) - \nu \sum_{j=1}^n P_j\r) \r\}\\
        &= \min_{\lambda \ge 0,  \nu}\l\{ \lambda \rgf{\rho} + \nu + \sum_{j=1}^n
        \max_{P_j \ge 0} \l(P_j \sum_{\tau=1}^L a_\tau  R^{i^j}_\tau - \lambda  Q_j
        \phi\l(\frac{P_j}{Q_j}\r)  - \nu P_j \r)\r\}\\
        &= \min_{\lambda \ge 0,  \nu}\l\{ \lambda \rgf{\rho} + \nu + \sum_{j=1}^n
        \max_{P_j \ge 0} \l(P_j \l(\sum_{\tau=1}^L a_\tau R^{i^j}_\tau - \nu \r) -
        \lambda  Q_j \phi\l(\frac{P_j}{Q_j}\r) \r)\r\} \\
        &= \min_{\lambda \ge 0,  \nu}\l\{ \lambda \rgf{\rho} + \nu + \lambda
        \sum_{j=1}^n Q_j \max_{t_j \ge 0}\l(t_j  \frac{R^{i^j}_\tau - \nu}{\lambda}
        - \phi\l(t_j\r) \r)\r\}\label{eq:sum_swap}\\
        &=\min_{\lambda \ge 0,  \nu}\l\{ \lambda \rgf{\rho} + \nu + \lambda
        \sum_{j=1}^n Q_j \max_{t_j \ge 0}\l(t_j s_j - \phi\l(t_j\r) \r)\r\}\\
        &= \min_{\lambda \ge 0,  \nu}\l\{ \lambda \rgf{\rho} + \nu + \lambda
        \sum_{j=1}^n Q_j \phi^*(s_j)\r\},
    \end{align}
    where $s_j = \frac{\sum_{\tau=1}^L a_\tau R^{i^j}_\tau - \nu}{\lambda}$ and
    $t_j = \frac{P_j}{Q_j}$. Note that we can replace the sum of maxima with a
    maximum of sums in~\eqref{eq:sum_swap} because the objective is separable
    over $j$. Finally, we require the dual feasibility
    constraint~\eqref{eq:dual_feasibility}\rgf{:}
    \begin{equation}\label{eq:dual_feasibility}
        s_j \le \l(\lim_{t \to \infty}\frac{\phi(t)}{t}\r) \fa j = 1,\dots,n\rgf{.}
    \end{equation}
    \rgf{This} ensures that $\phi^*$ does not grow to infinity. Consider
    $\phi^*(s_j) =\sup_{t \ge 0}\{s_j t - \phi(t)\}$. If $\frac{\phi(t)}{t} \to
    \infty$ as $t \to \infty$ then this constraint can be removed. If not, i.e.\
    $\lim_{t \to \infty}\frac{\phi(t)}{t} = \bar{s} < \infty$, then for $s >
    \bar{s}$ we have $\phi^*(s) = \infty$. Note that, according to the
    definition given by~\cite{BenTalHertog}, we have $0\phi^*(s/0) :=
    (0\phi^*)(s)$, which is zero if $s \le 0$ and $+\infty$ if $s > 0$.
    Therefore, combining with the outer problem, we have:
    \begin{align}
        \min_{y, R, \lambda, \nu} \ &\l\{ \lambda \rgf{\rho} + \nu + \lambda
        \sum_{j=1}^n Q_j \phi^*(s_j)\r\},\label{eq:reform_np_obj}\\
        \text{s.t. } & (2)-(8)\\
        &\lambda \ge 0 \\
        &\sum_{\tau=1}^L a_\tau R^{i^j}_\tau - \nu \le \lambda\l(\lim_{t \to
        \infty}\frac{\phi(t)}{t}\r) \fa j = 1,\dots,n.
    \end{align}
    
    \subsection{\texorpdfstring{Modified $\chi^2$-divergence}{}}
     
    Recall equation~(13), which states that for a modified $\chi^2$ divergence,
    we have: 
    $$d_{m\chi^2}(P, Q) = \sum_{j=1}^n \frac{(P_j - Q_j)^2}{Q_j}.$$
     
    \subsubsection{Reformulation}
    
    The conjugate of $\phi_{m\chi^2}$ is given by:
    \begin{align*}
        \phi^*_{m\chi^2}(s) &= \begin{cases}
            -1 &\text{ if } s < -2 \\
            s + \frac{s^2}{4} &\text{ if } s \ge -2
        \end{cases}\label{chi_conj} \\
        &= \max\l\{\frac{s_j}{2} + 1, 0\r\}^2 - 1.
    \end{align*}
    Using $\phi^*$ to represent $\phi^*_{m\chi^2}$ for shorthand, we can expand
    $\phi^*(s_j)$ in order to write:
    \begin{align}
        \lambda\phi^*(s_j) &= \lambda \l(\max\l\{\frac{s_j}{2} + 1, 0\r\}^2 - 1\r)\\
        &= \lambda \max\l\{\frac{\sum_{\tau=1}^L a_\tau R^{i^j}_\tau -
        \nu}{2\lambda} + 1, 0\r\}^2 - \lambda \\
        &= \frac{1}{4\lambda}\max\l\{\sum_{\tau=1}^L a_\tau R^{i^j}_\tau - \nu +
        2\lambda, 0\r\}^2 - \lambda.
    \end{align}
    In order to define $\phi^*(s_j)$ using convex quadratic constraints, we
    first need to remove the max operator from this expression. Hence, we define
    a dummy variable $z_j$ to represent the value of $\max\l\{\sum_{\tau=1}^L
    a_\tau R^{i^j}_\tau - \nu + 2\lambda, 0\r\}$. We enforce $z$'s value
    via~\eqref{eq:z_j_last1} and~\eqref{eq:z_j_last2}.
    \begin{align}
        z_j &\ge \sum_{\tau=1}^L a_\tau R^{i^j}_\tau - \nu + 2\lambda \fa j =
        1,\dots,n \label{eq:z_j_last1}\\
        z_j &\ge 0 \fa j = 1,\dots,n \label{eq:z_j_last2}.
    \end{align}
    Then, we can define another dummy variable $u_j = \frac{z_j^2}{\lambda} =
    4\lambda\phi^*(s_j) +  \lambda$. We enforce the value of $u_j$ using a conic
    quadratic constraint as follows:
    \begin{align}
        u_j &\ge \frac{z^2_j}{\lambda} \\
        \lambda u_j &\ge z^2_j \\
        (\lambda + u_j)^2 - (\lambda - u_j)^2 &\ge 4z^2_j,\\
        \sqrt{4z^2_j + (\lambda - u_j)^2} &\le (\lambda + u_j). 
    \end{align}
    Hence, with dummy variables $z_j, u_j$ for $j = 1,\dots,n$, we can
    reformulate our inner problem as:
    \begin{align}
        \min_{\lambda \ge 0,  \nu, z, u} \ &\l\{ \lambda (\rgf{\rho} - 1) + \nu +
        \frac{1}{4}\sum_{j=1}^n Q_j u_j \r\} \\
        &\sqrt{4z^2_j + (\lambda - u_j)^2} \le (\lambda + u_j) \fa j = 1, \dots, n
        \label{eq:NP_reform1}\\
        &z_j \ge \sum_{\tau=1}^L a_\tau R^{i^j}_\tau - \nu + 2\lambda \fa j =
        1,\dots,n\\
        &z_j \ge 0 \fa j = 1,\dots, n.\\
        &\lambda \ge 0\label{eq:NP_reform4}.
    \end{align}
    Therefore, combining with the \rgf{outer} model, we have:
    \begin{align}
        \min_{y, R, \lambda, \nu, z, u} &\l\{ \lambda (\rgf{\rho} - 1) + \nu +
        \frac{1}{4}\sum_{j=1}^n Q_j u_j \r\},\label{eq:CQP_obj}\\
        \text{s.t. } &(2)-(8),\\
        &\eqref{eq:NP_reform1}-\eqref{eq:NP_reform4}. \label{eq:CQP_last}
    \end{align}
    Note that, in the objective function, the $-\lambda$ comes from the fact
    that $\lambda \phi^*(s_j) = \frac{1}{4}u_j - \lambda$.
    
    \subsubsection{Extracting Worst-case Distribution}
    
    In order to find the worst-case distribution, we must extract it from the
    optimal values of $\lambda, \nu, z, u$. Denote the optimal solution of
    model~\eqref{eq:CQP_obj}-\eqref{eq:CQP_last} by $(y^*, R^*, \lambda^*,
    \nu^*, z^*, u^*)$. As discussed by~\cite{Love15}, the worst-case
    distribution $P^*$ satisfies:
    \begin{equation}\label{eq:optimality_conditions}
        \frac{P^*_j}{Q_j} \in \partial \phi^*(s^*_j), \qquad \sum_{j=1}^n Q_j
        \phi\l(\frac{P^*_j}{Q_j}\r) \le \rgf{\rho}, \qquad \sum_{j=1}^n P^*_j = 1.
    \end{equation}
    Here, the notation $\partial f(x)$ is the set of \textit{subgradients} of
    $f$ at $x$. Suppose that $\lambda^* > 0$ so that $s^*_j$ is defined.
    By~\cite{Love15}, if $\phi^*$ is differentiable then $(\phi^*)'(s^*_j)$ is a
    subgradient. This is true in our case, with $(\phi^*)'(s) =
    \max\l\{1+\frac{s}{2}, 0\r\}$. This derivative is non-negative, and hence
    always gives a feasible solution for $P^*_j$ by taking $P^*_j = Q_j
    (\phi^*)'(s^*_j)$ when $\lambda^* > 0$. In our experiments we only ever
    observed $\lambda^* > 0$ and hence $\phi^*(s^*_j)$ always gave a solution.
    For more detail on how to extract the solution when $\lambda^* = 0$,
    see~\cite{Love15}.
    
    \section{Further Analysis of Results}
    
    \subsection{The Effect of Workstacks on
    Solutions}\label{sec:effect_workstacks}
    
    In our experiments, we used only one value of the capacity $c$ but varied
    the workstacks $D$ to give a variety of possibilities for pulling forward.
    This was based on the number of pairs between which pulling forward was
    possible, i.e.\ $|\m{F}^+(c, D)|$ from Section~4.2. We give some examples of
    the values of $c - D$ and the corresponding $|\m{F}^+(c, D)|$ in
    Table~\ref{tab:num_pairs}.
    
    \begin{table}[htbp!]
        \centering
        \input{new_images_tables/fast_fp_examples}
        \caption{Examples of $c - D$ values and corresponding number of pairs}
        \label{tab:num_pairs}
    \end{table}
    
    \vspace{-1em}
    Any more pairs than 7 is not possible for $L = 5$ and $K = 2$. We present a
    summary of the results broken down by $|\m{F}^+(c, D)|$ in
    Table~\ref{tab:results_by_numpairs}. This table shows three quantities: the
    average time taken by each algorithm, the average gaps and the average
    number of pairs of days which had a positive pulling forward decision. The
    table shows that we did not have any more non-zero decisions than 1, from
    any algorithm, until $|\m{F}^+(c, D)|$ reached its maximum value of 7. Days
    1 and 2 \rgf{were} typically prioritised for rollover reduction via pulling
    forward in these cases. This is because jobs due on these days have the
    potential to roll over the most times. However, when $|\m{F}^+(c, D)| \ =
    7$, we see between 1 and 6 pairs of days having a non-zero pulling forward
    decision.
    
    \begin{table}[htbp!]
        \centering
        \small
        \input{new_images_tables/final2_results_by_fp}
        \caption{Results by $|\m{F}^+(c, D)|$.}
        \label{tab:results_by_numpairs}
    \end{table}
    
    The APGs are also shown in Table~\ref{tab:results_by_numpairs}. From this,
    we can see a number of results. Firstly, we see that the average time taken
    by each algorithm apart from AO \rgf{was} increasing in $|\m{F}^+(c, D)|$.
    This can be expected, since more feasible pairs leads to a more complex
    feasible region. Furthermore, AO perform\rgf{ed the} worse in selecting $y$
    as $|\m{F}^+(c, D)|$ increases. This is likely because reducing the set of
    intakes leads to less accurate estimates of the expected rollover.
    Interestingly, CS \rgf{did} not suffer from the same issue. In fact, for
    $|\m{F}^+(c, D)| \ = 7$, CS ha\rgf{d} an average $y$-APG of 0.0\% and for
    all values of $|\m{F}^+(c, D)|$ this value \rgf{was} below 0.051\%. This is
    because CS does not employ dimension reduction to the set of intakes like AO
    does. As might be expected, there is no clear pattern in the $p$-APGs. For
    AO and CS, this value \rgf{was} highest when $|\m{F}^+(c, D)| \ = 5$ and
    lowest when $|\m{F}^+(c, D)| \ = 7$. Finally, the final column shows the
    average and maximum numbers of pairs $(\tau_1, \tau_2)$ that had $y_{\tau_1,
    \tau_2} > 0$ under each algorithm. The results for $|\m{F}^+(c, D)| \ = 7$
    suggest that NP's solution \rgf{was} slightly less conservative than P's
    solution on average. We study this in more detail in Section~5.5.
    Interestingly, NP \rgf{took} almost as long as P in these instances. CS\_opt
    again ha\rgf{d} all zero gaps and APGs, and its times taken \rgf{were} no
    more affected by $|\m{F}^+(c, D)|$ than the times taken by CS.
    
    \subsection{Comparison with Robust Optimisation Solutions} \label{sec:RO}
    
    In this section, we compare the DRO decisions and objectives with those
    resulting from the robust optimisation (RO) version of the model. The RO
    model is obtained by replacing the inner objective with the maximisation of
    the total rollover cost over all intake vectors. The first result that we
    find is that the intake vector responsible for the worst-case cost for the
    chosen $y$ value was always $i^{\max}$. This shows that the RO model can be
    solved simply by assuming that all intakes take their maximum values at all
    times. As well as this, the RO model pulled forward less than the DRO model
    in 227 (82\%) of our 279 instances. The RO solution also had a higher cost
    than the DRO solution in 269 (97\%) of instances. This can be expected due
    to the way that their \rgf{objective functions} differ. These two facts
    \rgf{support} our claim that the RO model is more conservative than the DRO
    model.
    
    We present some more detailed results in Table~\ref{tab:DRO_RO_results}.
    This table compares the objective values, pulling forward decisions and
    times taken from the three models. Firstly, note that RO \rgf{took} around
    16 seconds on average. RO also pull\rgf{ed} forward less than DRO.
    Specifically, it pulled forward 1.3 jobs less than DRO, on average. Also,
    DRO pulled forward a maximum of 8 jobs whereas RO only pulled forward a
    maximum of 7 jobs. Furthermore, the objective values from RO were
    significantly higher than DRO. Comparing the RO objective with the DRO
    objective, we see that RO's objective values were around 9.5 higher than
    DRO's on average. This corresponds to almost a 200\% increase in objective
    value. The $y$-gap and $y$-APGs assess the expected costs from RO's
    decisions when evaluated by DRO's objective function. This suggests that
    RO's decisions would result in around 2 more jobs being expected to roll
    over in the worst case than DRO's solution. 
    
    \begin{table}[htbp!]
        \centering
        \input{new_images_tables/RO_summary}
        \caption{Comparison of results from RO model with DRO solutions}
        \label{tab:DRO_RO_results}
    \end{table}
    
    \vspace{-1em}
    As already noted, RO is equivalent to the deterministic model under the
    assumption that $I = i^{\max}$ with probability 1. The results from this
    model are shown in the ``RO det.''\ column. This shows that this model took
    0.01 seconds to build and solve, on average. Hence, our results indicate
    that the inclusion of the rollover constraints for the RO model \rgf{led} to
    around a 16 second increase in solution times. The inclusion of the expected
    value constraints for the DRO model result\rgf{ed} in over 1 minute of
    additional solution time. Table~\ref{tab:DRO_RO_results} also shows that RO
    had an objective value that was three times larger than DRO's, on average.  
    
    From the results presented here, we can conclude three main results.
    Firstly, RO is more conservative than DRO for this problem, since it
    pull\rgf{ed} forward fewer jobs on average. Secondly, RO results in
    significantly higher costs for the same $y$ decision. However, the third
    conclusion is that RO is much faster than DRO. This indicates that the main
    factor affecting solution times for DRO is the inclusion of the expected
    value constraints.
    
    % \subsection{Robustness of DRO Solutions}
    
    % In this section, we summarise the robustness of our model's solutions
    % against changes in the demand distribution, and compare them with the
    % solutions that would be obtained if the distribution were known. Since we
    % did not assume a true distribution earlier in the paper, in order to
    % assess robustness, we have carried out some additional experiments.
    % Specifically, we solved the same instances as before, but with some slight
    % changes in methodology. This time, for each instance we randomly generated
    % a true parameter $p^0$ and used this to generate the MLE $\hat{p}$ and the
    % corresponding confidence set $\Theta_{\alpha}$. In addition, we solved the
    % stochastic program that results from the true distribution in order to
    % obtain the optimal cost under this distribution. This allows us to assess
    % the loss in performance associated with the DRO solution. Furthermore, we
    % evaluated the cost of the DRO solution under every distribution in the
    % ambiguity set, in order to see how much it varies when the distribution
    % changes. We present the results here, due to space considerations.

    \section{A Benders Decomposition Approach}
    
    Our CS\_opt algorithm can be viewed as a specialised Benders
    decomposition~\citep{benders} approach that solves the distribution
    separation problem as a residual problem. However, it does not require us to
    create the dual of the distribution separation problem, and in our case we
    can simply solve this problem by enumeration. For comparison, we now present
    a classical Benders decomposition approach in order to explain why CS\_opt
    is preferred.
    
    \subsection{Residual Problem and its Dual}
    
    We create the Benders residual problem by taking $y$ as \rgf{the} master
    problem variable and $R, t$ as \rgf{the} subproblem variables. This is
    because the model's complexity comes from $R$ and $t$, not $y$. For a fixed
    $y = \bar{y}$, the residual problem can be written as:
    \begin{align}
        \min_{R, t} \quad & t,\\
        \text{s.t. } &R^i_1 \ge i_1 + \sum_{\tau_1 = 2}^{\min\{1 + K, L\}}
        y_{\tau_1, 1} - \l(c_1 - D_1\r) \fa i \in \m{I}\\
        &R^i_\tau - R^i_{\tau - 1} \ge i_\tau + \sum_{\tau_1 = \tau + 1}^{\min\{\tau
        + K, L\}} y_{\tau_1, \tau} - \l(c_\tau - D_\tau + \sum_{\tau_2 = \max\{\tau
        - K, 1\}}^{\tau -1} y_{\tau, \tau_2}\r) \nonumber\\
        & \fa \tau = 2, \dots,L-1 \fa i \in \m{I},\\
        &R^i_L - R^i_{L - 1} \ge i_L - \l(c_L - D_L + \sum_{\tau_2 = \max\{L - K,
        1\}}^{L - 1} y_{\tau, \tau_2}\r) \fa i \in \m{I},\\
        & t - \sum_{\tau=1}^L a_\tau \sum_{i \in \m{I}} \P(I = i| p) R^i_\tau \ge 0
        \fa p \in \Theta.
    \end{align}
    This model has $m = L|\m{I}| + |\Theta|$ constraints. Hence, we have dual
    variables $u_{j, \tau}$ for $j = 1,\dots, |\m{I}|$ and $\tau = 1,\dots, L$,
    and $v_k$ for $k = 1,\dots, |\Theta|$. The model has $L|\m{I}| + 1$
    variables, and so we have $L|\m{I}|+1$ constraints in the dual. The dual is
    given by:
    \begin{align}
        \max_{x} &\sum_{\tau=1}^L\sum_{j = 1}^n b_{j, \tau}(\bar{y}) u_{j, \tau} +
        \sum_{k=1}^{|\Theta|} \tilde{b}_k(\bar{y}) v_k \label{eq:subproblem_1}\\
        \text{s.t. } R^i_1: & \ \ u_{j, 1} - u_{j, 2} -
        \sum_{k=1}^{|\Theta|}a_1\P(I=i^j| p^k) v_k \le 0 \fa j = 1,\dots,|\m{I}|,\\
        R^i_\tau: & \ \ u_{j, \tau} - u_{j, \tau+1} -
        \sum_{k=1}^{|\Theta|}a_\tau\P(I=i^j| p^k)v_k \le 0 \fa j = 1,\dots,|\m{I}|
        \fa \tau = 2,\dots, L-1,\\
        R^i_L: & \ \ u_{j, L} - \sum_{k=1}^{|\Theta|}a_L\P(I=i^j| p^k)v_k \le 0 \fa
        j = 1,\dots,|\m{I}|,\\
        t: & \ \ \sum_{k=1}^{|\Theta|} v_k \le 1,\label{eq:subproblem_last}
    \end{align}
    where $b_{j, \tau}(\bar{y})$ and $\tilde{b}_k(\bar{y})$ are defined as:
    \begin{align}
        b_{j, 1}(\bar{y}) &=  i_1 + \sum_{\tau_1 = 2}^{\min\{1 + K, L\}}
        \bar{y}_{\tau_1, 1} - \l(c_1 - D_1\r) \fa i \in \m{I}\\
        b_{j, \tau}(\bar{y}) &=  i_\tau + \sum_{\tau_1 = \tau + 1}^{\min\{\tau + K,
        L\}} \bar{y}_{\tau_1, \tau} - \l(c_\tau - D_\tau + \sum_{\tau_2 = \max\{\tau
        - K, 1\}}^{\tau -1} \bar{y}_{\tau, \tau_2}\r) & \fa \tau = 2, \dots,L-1 \fa
        i \in \m{I},\\
        b_{j, L}(\bar{y}) &=  i_L - \l(c_L - D_L + \sum_{\tau_2 = \max\{L - K,
        1\}}^{L - 1} \bar{y}_{\tau, \tau_2}\r) \fa i \in \m{I},\\
        \tilde{b}_k(\bar{y}) &= 0 \fa k = 1,\dots,|\Theta|.
    \end{align}
    
    \subsection{Benders Decomposition Algorithm}
    
    Our Benders decompositon algorithm is as follows.
    \begin{enumerate}
        \item Initialise $\varepsilon$, $LB = - \infty$, $UB = \infty$. Set
        feasible region for $z$ as $Z =  \mathbb{R}^+$. Set feasible region for
        $y$ as $Y$, where $y \in Y$ indicates that $y$ is feasible for the model
        in (2)-(8).
        \item While $UB - LB > \varepsilon$:
        \begin{enumerate}
            \item Solve master problem:
            \begin{equation}
                \min_{z \in Z, y \in Y} z \label{eq:master}
            \end{equation}
            to get a solution $\bar{y}$ and objective value $z^{\text{M}}$.
            \item Set $LB = z^{\text{M}}$.
            \item Solve Benders
            subproblem~\eqref{eq:subproblem_1}-\eqref{eq:subproblem_last} with
            $y = \bar{y}$ to get a solution $\bar{u}$, $\bar{v}$ with objective
            $z^{\text{S}}$.
            \item If subproblem is unbounded, add feasibility cut:
            $$ \sum_{\tau=1}^L\sum_{j = 1}^n b_{j, \tau}({y}) u_{j, \tau} +
            \sum_{k=1}^{|\Theta|} \tilde{b}_k({y}) v_k \le 0$$ to $Y$.
            \item If subproblem is optimal, add optimality cut:
            $$z \ge \sum_{\tau=1}^L\sum_{j = 1}^n b_{j, \tau}({y}) u_{j, \tau} +
            \sum_{k=1}^{|\Theta|} \tilde{b}_k({y}) v_k$$ to $Z$.
            \item If $z^{\text{S}} < UB$ then set $UB = z^{\text{S}}$.
        \end{enumerate}
        \item Find index of binding $t$ constraint from the subproblem and use
        this to find worst-case $p$.
        \item Return $y$, $p$.
    \end{enumerate}
    In the following section we will show that this approach is slow compared
    with CS\_opt.
    
    \subsection{Results}
    
    We tested the Benders algorithm on each of our 279 instances, for
    $\varepsilon \in \{0.01, 10^{-6}, 10^{-8}\}$. We present the results in
    Table~\ref{tab:benders}. From these results, it is clear that $\varepsilon =
    10^{-8}$ was required for $y$ optimality. However, with this $\varepsilon$,
    the Benders algorithm took almost 6 minutes to solve, on average. In one
    instance, the algorithm timed out as it took longer than 4 hours. In
    comparison with CS\_opt, which takes approximately 17 seconds on average,
    this version of Benders decomposition is very slow.
    \begin{table}[htbp!]
        \centering
        \input{new_images_tables/benders_results}
        \caption{Results of Benders algorithm}
        \label{tab:benders}
    \end{table}

    \newpage

    \section{Large Results Tables}\label{sec:tables_appendix}
    
    \subsection{Results by \texorpdfstring{$|\Theta|$}{}}\label{sec:U_table}
    \begin{table}[htbp!]
        \centering
        \footnotesize
        \input{new_images_tables/final2_perf_by_Theta}
        \caption{Summary of results and times taken by $N$ and
        $n_{\text{probs}}$. Referred to in Section~5.3.}
        \label{tab:results_by_AS_size}
    \end{table}
    
    \newpage
    
    \subsection{Results by \texorpdfstring{$|\m{I}|$}{}}\label{sec:I_table}
    
    \begin{table}[htbp!]
        \centering
        \footnotesize
        \input{new_images_tables/final2_results_by_US_size}
        \caption{Summary of results and times taken by size of $\m{I}$. Referred
        to in Section~5.3.}
        \label{tab:results_by_US_size}
    \end{table}
    
    \vspace{1em}
    
    \newpage
    
    \section{Tables of Notation}\label{sec:notation_tables}
    
    \subsection{General Model Notation}
    
    \begin{table}[htbp!]
        \footnotesize
        \begin{tabularx}{\columnwidth}{X|X}
            \hline
            
            \textbf{Notation} & \textbf{Meaning} \\
            \hline
            $L$ & Number of days in a plan \\
            $K$ & Maximum number of days a job can be pulled forward \\
            $\tau$, $\tau_1$, $\tau_2$ & A day in the plan, value in $\{1,\dots,
            L\}$ \\
            $y_{\tau_1, \tau_2}$ & Number of jobs to pull forward from day
            $\tau_1 \in \{2,\dots,L\}$ to $\tau_2 \in \{\max{\tau_1 - K, 1},
            \dots, \tau_1 - 1\rgf{\}}$. \\
            $R_\tau$ & Number of jobs to roll over from day $\tau$ to $\tau +
            1$.\\
            $a_\tau$ & Cost of a job rolling over from day $\tau$ to $\tau +
            1$.\\
            $c_\tau$ & Number of hours of capacity available on day $\tau$.\\
            $D_\tau$ & Number of jobs currently due on day $\tau$.\\
            $\mathbb{N}_0$ & Set of non-negative integers.\\
            $I_\tau$ & Random variable representing number of jobs arriving
            between the time of planning and day $\tau$ that will be due on day
            $\tau$ (intake).\\
            $i_\tau$ & Realisation of $I_\tau$.\\
            $R^i$ & Realisation of $R = (R_1, \dots, R_L)$ corresponding to
            realisation $i$ of $I$.\\
            $\m{I}_\tau$ & Set of all possible realisations of $I_\tau$.\\
            $\m{I}$ & Set of all possible realisations of the vector of intakes
            $I$.\\
            $\m{P}$ & General ambiguity set constaining distributions of
            intake.\\
            $P$ & A discrete probability distribution over the set of intakes
            $\m{I}$. \\
            $Q$ & Nominal distribution of intake.\\
            $i^{\max}_\tau$ & The maximum value $I_\tau$ can take.\\
            $p_\tau$ & A variable representing success probability parameter of
            the binomial distribution of intake $I_\tau$.\\
            $p^0_\tau$ & True success probability of intake $I_\tau$.\\
            $\hat{p}_\tau$ & MLE of $p^0_\tau$ taken from $N$ samples of
            $I_\tau$.\\
            $P^p$ & Binomial distribution of intake with success probability
            parameter $p$.\\
            $\hat{P}$ & MLE of distribution resulting from $p = \hat{p}$.\\
            $\m{P}_{\Theta}$ & Set of all probability distributions $P$ that are
            binomial with a value $p \in \Theta$. \\
            $\Theta$ & Set of vectors $p$ obtained from a distribution $P$ in
            $\m{P}_{\Theta}$.\\
            $\Theta_\alpha$ & $100(1-\alpha)\%$ confidence set for $p^0$ around
            the MLE $\hat{P}$. 
        \end{tabularx}
        \label{tab:general_notation}
        \caption{General model notation from Section~3.}
    \end{table}
    
    \newpage 
    
    \subsection{Non-parametric Model Notation}
    
    \begin{table}[htbp!]
    \begin{tabularx}{\columnwidth}{X|X}
        \hline
        \textbf{Notation} & \textbf{Meaning} \\
        \hline
        $\phi$ & $\phi$-divergence function.\\
        $d_\phi$ & $\phi$-divergence measure resulting from $\phi$-divergence
        function $\phi$.\\
        $\phi^*$ & Conjugate of $\phi$-divergence function $\phi$.\\
        $\phi_{m\chi^2}$ & $\phi$-divergence function for modified $\chi^2$
        distance.\\
        $\chi^2_{k, 1-\alpha}$ & $100(1-\alpha)\%$ percentile of $\chi^2$
        distribution with $k$ degrees of freedom.\\
        $\m{P}_\rho$ & Non-parametric confidence set for true distribution
        $P^0$.\\
        $\lambda, \nu$ & Lagrange multipliers for SQP reformulation of NP
        model.\\
        $u_j, z_j, s_j$ & Dummy variables used to reformulate NP model.\\
        $\rgf{\rho}$ & Maximum distance, measured by $d_\phi$, from $Q$ that $P$
        can be under the non-parametric model.\\
        $\partial{f(x)}$ & Set of subgradients of a function $f$ at a point
        $x$.\\
        $\cdot^*$ & Optimal value of $\cdot$ under the non-parametric
        model\rgf{, for $\cdot \in \{R, s, P, y, \nu, \lambda\}$}.
        \end{tabularx}
    \label{tab:NP_notation}
    \caption{Notation used in the non-parametric model in Section~3.3}
    \end{table}
    
    \newpage
    \subsection{CS/CS\_opt/AO Notation}
    
    \begin{table}[htbp!]
        \begin{tabularx}{\columnwidth}{X|X}
            \hline
            \textbf{Notation} & \textbf{Meaning} \\
            \hline
            $k$ (repeated dummy variable) & Index for the iteration of
            CS/CS\_opt algorithm that we are currently carrying out.\\
            $p^{\max}_\tau$ & Maximum value that $p_\tau$ takes over $p \in
            \Theta$.\\
            $\Theta^{\max}_\tau$ & Set of $p$ parameters such that $p_\tau$ is
            maximised.\\
            $\Theta^{\text{ext}}$ & Set of extreme \rgf{parameters} used by
            CS.\\
            $\tilde{\Theta}$ & General ambiguity set used by CS algorithms.
            $\tilde{\Theta} = \Theta$ for CS\_opt and $\tilde{\Theta} =
            \Theta^{\text{ext}}$ for CS.\\
            $k^{\max}$ & Maximum number of iterations of CS/CS\_opt algorithm
            allowed to run. \\
            ${\Theta}^{k}$ & Current subset of ${\Theta}$ being used at
            iteration $k$ of CS/CS\_opt. \\
            $y^k$ & Pulling forward decision generated by solving outer problem
            at iteration $k$ of CS/CS\_opt.\\
            $p^k$ & Probability vector generated by solving distribution
            separation problem at iteration $k$ of CS/CS\_opt.\\
            $\varepsilon$ & Optimality tolerance of CS/CS\_opt algorithm.\\
            $t^k$ & Objective value of problem obtained by solving outer problem
            at iteration $k$ of CS.\\
            $\beta$ & Minimum probability an intake must have of occurring in
            order to be used in the AO algorithm. \\
            $\tilde{\m{I}}$ & Set of intakes with probability of occurring
            higher than $\beta$.\\
        \end{tabularx}
    \label{tab:CS_AO_notation}
    \caption{Notation used in CS/AO Algorithms (Section~3.6)}
    \end{table}
    
    \newpage
    
    \subsection{Input Parameter and Results Notation}
    \begin{table}[htbp!]
        \begin{tabularx}{\columnwidth}{X|X}
            \hline
            \textbf{Notation} & \textbf{Meaning} \\
            \hline
            $\mathcal{F}$ & Set of pairs of days between which pulling forward
            is allowed.\\
            $\mathcal{F}^+(c, D)$ & Set of pairs of days between which pulling
            forward is feasible given $c$ and $D$.\\
            $n(i^{\max})$ & Number of days with maximum intake higher than
            remaining capacity given $i^{\max}$, $c$, and $D$.\\
            $n_{\text{probs}}$ & Number of values each probability in $p$ can
            take in our discretised ambiguity set. \\
            $\hat{i}$ & MLE of mean intake vector.\\
            $\rho$ & Maximum distance from $\hat{P}$ we allow $P$ to be under
            NP, measured by the chosen $\phi$-divergence.\\
      \end{tabularx}
    \label{tab:Input_notation}
    \caption{Input parameter notation used in Section~4}
    \end{table}
    
    \begin{table}[htbp!]
        \begin{tabularx}{\columnwidth}{X|X}
            \hline
            \textbf{Notation} & \textbf{Meaning} \\
            \hline
            $f(y, p)$ & Shorthand for expected rollover cost given pulling
            forward decision $y$ and distribution parameter $p$.\\
            $x$ & An algorithm, namely in $\{\text{S\&S}, \text{CS},
            \text{AO}\}$.\\
            $y^x, p^x$ & $y, p$ solution obtained by algorithm $x$.\\
            $g_p(y^x, p^x)$ & $p$-gap of algorithm $x$'s solution. The
            difference between the worst-case expected cost for $y^x$ and the
            expected cost obtained by the algorithm.\\
            $z^*$ & Overall optimal objective value.\\
            $g_p(y^x)$ & $y$-gap. Difference between worst-case expected cost
            for $y^x$ over all distributions and the overall optimal objective
            value.
        \end{tabularx}
    \label{tab:results_notation}
    \caption{Results analysis notation from Section~5}
    \end{table}
    
    \end{appendices}

%% file: new_images_tables/fast_fp_examples.tex
\begin{tabular}{lr}
\toprule
  $c - D$ &  $|\mathcal{F}^+(c, D)|$ \\
\midrule
(8, -15, -15, 8, -15) &                        3 \\
(8, -15, 8, 8, 8) &                        5 \\
(8, 8, 8, 8, 8) &                        7 \\
\bottomrule
\end{tabular}

%% file: new_images_tables/final2_results_by_fp.tex
\begin{tabular}{lllllll}
\toprule
  &     &    & Avg. $p$-APG & Avg. $y$-APG &   Avg. t.t. & (Avg., Max) Non-zeros \\
$|\mathcal{F}^+(c, D)|$ & Count & Algorithm &              &              &             &                       \\
\midrule
3 & 189 & opt &        0.0\% &        0.0\% &  0:01:00.38 &              (1.0, 1) \\
  %&     & SnS &     0.0398\% &     1.5159\% &  0:00:06.29 &              (1.0, 1) \\
  &     & CS &     0.0415\% &     0.0029\% &  0:00:07.58 &              (1.0, 1) \\
  &     & CS_opt &        0.0\% &        0.0\% &  0:00:09.75 &              (1.0, 1) \\
  &     & AO &     0.0084\% &     0.0059\% &  0:00:24.74 &              (1.0, 1) \\
  &     & NP &            - &            - &  0:00:02.62 &              (1.0, 1) \\
  \midrule
5 & 45  & opt &        0.0\% &        0.0\% &  0:02:12.72 &              (1.0, 1) \\
  %&     & SnS &     0.2936\% &     6.8772\% &  0:00:10.47 &              (1.0, 1) \\
  &     & CS &     0.3465\% &     0.0506\% &  0:00:25.62 &              (1.0, 1) \\
  &     & CS_opt &        0.0\% &        0.0\% &  0:00:29.35 &              (1.0, 1) \\
  &     & AO &     0.0052\% &     0.2903\% &  0:00:19.04 &              (1.0, 1) \\
  &     & NP &            - &            - &   0:00:32.2 &              (1.0, 1) \\
  \midrule
7 & 45  & opt &        0.0\% &        0.0\% &  0:02:44.25 &              (2.0, 4) \\
  %&     & SnS &        0.0\% &   115.1682\% &  0:00:26.84 &              (0.8, 4) \\
  &     & CS &        0.0\% &        0.0\% &  0:01:27.85 &              (2.0, 4) \\
  &     & CS_opt &        0.0\% &        0.0\% &  0:01:31.68 &              (2.0, 4) \\
  &     & AO &        0.0\% &     0.5332\% &  0:00:21.16 &              (2.0, 4) \\
  &     & NP &            - &            - &  0:01:53.93 &              (2.8, 6) \\
\bottomrule
\end{tabular}

%% file: new_images_tables/RO_summary.tex
\begin{tabular}{llll}
\toprule
{} &     RO Det. &          RO &         DRO \\
\midrule
Avg. Obj. Gap                                    &       9.499 &       9.499 &           0 \\
Avg. \% Obj. Gap                                 &   199.662\% &   199.662\% &         0\% \\
Avg. $y$-gap                                     &       1.851 &       1.851 &           0 \\
Avg. $y$-APG                                     &     2.616\% &     2.616\% &         0\% \\
Avg. $\sum_{\tau_1, \tau_2}  y_{\tau_1, \tau_2}$ &           4 &           4 &       5.308 \\
Max. $\sum_{\tau_1, \tau_2}  y_{\tau_1, \tau_2}$ &           7 &           7 &           8 \\
Avg. t.t.                                        &  0:00:00.01 &  0:00:15.89 &  0:01:22.85 \\
\bottomrule
\end{tabular}

%% file: new_images_tables/benders_results.tex
\begin{tabular}{rrlrlll}
\toprule
 $\varepsilon$ &  Avg. $p$-gap & Avg. $p$-APG &  Avg. $y$-gap & Avg. $y$-APG &     Avg. t.t. &     Max t.t. \\
\midrule
 1e-08 &           0.0 &        0.0\% &       -0.0000 &        0.0\% &   0:05:57.798 &  4:00:05.011 \\
 1e-06 &           0.0 &        0.0\% &        0.0029 &        0.0\% &  0:04:42.6717 &  1:24:55.625 \\
 1e-02 &           0.0 &        0.0\% &        0.0141 &    10.1984\% &  0:04:36.6348 &  1:21:58.415 \\
\bottomrule
\end{tabular}

%% file: new_images_tables/final2_perf_by_Theta.tex
\begin{tabular}{llll|lll}
\toprule 
$(N, n_{\text{probs}})$ & Avg. $|\Theta|$ & Count & Algorithm &           Avg. $p$-APG & Avg. $y$-APG &   Avg. t.t. \\
\midrule
$(100, 5),(100,10)$ & 1.000    & 62 & P &        0.0\% &        0.0\% &  0:00:14.57 \\
         & &    & CS &        0.0\% &        0.0\% &     0:00:16 \\
         & &    & CS_opt &        0.0\% &        0.0\% &  0:00:15.81 \\
         & &    & AO &        0.0\% &      0.248\% &  0:00:01.32 \\
         & &    & NP &            - &            - &  0:00:25.58 \\
         \midrule
$(50, 5)$ &1.419    & 31 & P &        0.0\% &        0.0\% &  0:00:15.11 \\
         & &    & CS &        0.0\% &        0.0\% &  0:00:16.53 \\
         & &    & CS_opt &        0.0\% &        0.0\% &  0:00:16.46 \\
         & &    & AO &        0.0\% &      0.248\% &  0:00:01.35 \\
         & &    & NP &            - &            - &  0:00:27.13 \\
         \midrule
$(50, 10)$ & 14.419   & 31 & P &        0.0\% &        0.0\% &  0:00:14.23 \\
         & &    & CS &        0.0\% &        0.0\% &  0:00:15.05 \\
         & &    & CS_opt &        0.0\% &        0.0\% &  0:00:15.42 \\
         & &    & AO &        0.0\% &      0.248\% &   0:00:01.4 \\
         & &    & NP &            - &            - &  0:00:25.22 \\
         \midrule
$(10, 5)$ & 14.742   & 31 & P &        0.0\% &        0.0\% &  0:00:16.14 \\
         & &    & CS &        0.0\% &        0.0\% &  0:00:18.97 \\
         & &    & CS_opt &        0.0\% &        0.0\% &  0:00:19.42 \\
         & &    & AO &        0.0\% &     0.0004\% &  0:00:01.76 \\
         & &    & NP &            - &            - &  0:00:24.66 \\
         \midrule
$(100, 15)$ & 16.871   & 31 & P &        0.0\% &        0.0\% &  0:00:16.24 \\
         & &    & CS &        0.0\% &        0.0\% &  0:00:16.48 \\
         & &    & CS_opt &        0.0\% &        0.0\% &   0:00:16.5 \\
         & &    & AO &     0.0452\% &     0.2391\% &  0:00:01.57 \\
         & &    & NP &            - &            - &  0:00:26.12 \\
         \midrule
$(50, 15)$ & 93.129   & 31 & P &        0.0\% &        0.0\% &  0:00:24.51 \\
         & &    & CS &     0.5029\% &     0.0734\% &   0:00:20.8 \\
         & &    & CS_opt &        0.0\% &        0.0\% &   0:00:19.9 \\
         & &    & AO &     0.0105\% &        0.0\% &   0:00:02.9 \\
         & &    & NP &            - &            - &  0:00:27.16 \\
         \midrule
$(10, 10)$ & 504.226  & 31 & P &        0.0\% &        0.0\% &  0:00:59.36 \\
         & &    & CS &     0.0339\% &        0.0\% &  0:00:17.85 \\
         & &    & CS_opt &        0.0\% &        0.0\% &  0:00:19.76 \\
         & &    & AO &        0.0\% &        0.0\% &  0:00:12.77 \\
         & &    & NP &            - &            - &  0:00:23.47 \\
         \midrule
$(10, 15)$ & 4301.645 & 31 & P &        0.0\% &        0.0\% &  0:09:30.88 \\
         & &    & CS &      0.219\% &     0.0176\% &  0:00:19.66 \\
         & &    & CS_opt &        0.0\% &        0.0\% &  0:00:42.47 \\
         & &    & AO &     0.0031\% &        0.0\% &  0:02:11.27 \\
         & &    & NP &            - &            - &  0:00:23.22 \\
\bottomrule
\end{tabular}

%% file: new_images_tables/final2_results_by_US_size.tex
\begin{tabular}{lll|lll}
\toprule

$|\mathcal{I}|$ & Count & Algorithm & Avg. $p$-APG & Avg. $y$-APG &   Avg. t.t. \\
\midrule
392   & 27 & P &        0.0\% &        0.0\% &  0:00:33.63 \\
      &    & CS &     0.0234\% &        0.0\% &  0:00:06.37 \\
      &    & CS_opt &        0.0\% &        0.0\% &  0:00:08.01 \\
      &    & AO &      0.002\% &        0.0\% &  0:00:23.53 \\
      &    & NP &            - &            - &  0:00:00.84 \\
      \midrule
512   & 45 & P &        0.0\% &        0.0\% &  0:00:31.88 \\
      &    & CS &     0.0952\% &     0.0024\% &  0:00:06.54 \\
      &    & CS_opt &        0.0\% &        0.0\% &  0:00:08.05 \\
      &    & AO &        0.0\% &        0.0\% &  0:00:20.86 \\
      &    & NP &            - &            - &  0:00:00.78 \\
      \midrule
567   & 45 & P &        0.0\% &        0.0\% &  0:00:37.09 \\
      &    & CS &     0.0572\% &     0.0107\% &  0:00:06.61 \\
      &    & CS_opt &        0.0\% &        0.0\% &  0:00:08.28 \\
      &    & AO &     0.0116\% &        0.0\% &  0:00:25.75 \\
      &    & NP &            - &            - &  0:00:00.98 \\
      \midrule
2187  & 27 & P &        0.0\% &        0.0\% &   0:01:42.9 \\
      &    & CS &     0.0536\% &        0.0\% &  0:00:08.72 \\
      &    & CS_opt &        0.0\% &        0.0\% &  0:00:11.69 \\
      &    & AO &     0.0259\% &     0.0005\% &  0:00:25.93 \\
      &    & NP &            - &            - &  0:00:04.25 \\
      \midrule
2592  & 45 & P &        0.0\% &        0.0\% &     0:01:42 \\
      &    & CS &     0.0043\% &        0.0\% &  0:00:09.69 \\
      &    & CS_opt &        0.0\% &        0.0\% &  0:00:12.82 \\
      &    & AO &     0.0061\% &     0.0246\% &  0:00:26.58 \\
      &    & NP &            - &            - &  0:00:06.18 \\
      \midrule
8192  & 45 & P &        0.0\% &        0.0\% &  0:02:12.72 \\
      &    & CS &     0.3465\% &     0.0506\% &  0:00:25.62 \\
      &    & CS_opt &        0.0\% &        0.0\% &  0:00:29.35 \\
      &    & AO &     0.0052\% &     0.2903\% &  0:00:19.04 \\
      &    & NP &            - &            - &   0:00:32.2 \\
      \midrule
20000 & 45 & P &        0.0\% &        0.0\% &  0:02:44.25 \\
      &    & CS &        0.0\% &        0.0\% &  0:01:27.85 \\
      &    & CS_opt &        0.0\% &        0.0\% &  0:01:31.68 \\
      &    & AO &        0.0\% &     0.5332\% &  0:00:21.16 \\
      &    & NP &            - &            - &  0:01:53.93 \\
\bottomrule
\end{tabular}